\documentclass[12pt]{article}
\usepackage{amsmath}
\usepackage{amssymb}
\usepackage{amsthm}
\usepackage{graphicx}
\usepackage[usenames]{color}

\numberwithin{equation}{section}

\newcommand{\bfu}{{\boldsymbol u}}

\newcommand{\bfw}{{\boldsymbol w}}
\newcommand{\bfz}{{\boldsymbol z}}
\newcommand{\bff}{{\boldsymbol f}}

\newcommand{\bfr}{{\boldsymbol r}}
\newcommand{\bfa}{{\boldsymbol a}}

\newcommand{\bfG}{{\boldsymbol G}}

\newcommand{\bfY}{{\boldsymbol Y}}

\newtheorem{thm}{Theorem}[section]

\newtheorem{rem}[thm]{Remark}
          \newenvironment{req}{\begin{rem}\rm}{\end{rem}}

\newtheorem{prop}[thm]{Proposition} 
\newtheorem{lem}[thm]{Lemma} 
\newtheorem{cor}[thm]{Corollary}
\newtheorem{definition}[thm]{Definition} 
        \newenvironment{defn}{\begin{definition}\rm}{\end{definition}}
\newtheorem{example}[thm]{Example} 
        \newenvironment{ex}{\begin{example}\rm}{\end{example}}

\renewcommand{\qed}{\ifmmode$\Box$\else{\unskip\nobreak\hfil
	\penalty50\hskip1em\null\nobreak\hfil$\Box$
	\parfillskip=0pt\finalhyphendemerits=0\endgraf}\fi}

\begin{document}
\noindent

\vspace*{1cm}
\begin{center}
 {\Large \bf PARABOLIC QUASI-VARIATIONAL INEQUALITIES (I)\\[0.3cm]
     -- SEMIMONOTONE OPERATOR APPROACH --} \vspace{1cm}
\end{center}

\begin{center}
     {\sc Maria Gokieli}{\footnote{Corresponding author}}

Faculty of Mathematics and Natural Sciences. School
 of Exact Sciences \\
Cardinal Stefan Wyszynski University, Warsaw, Poland\\
 (E-mail: m.gokieli@uksw.edu.pl)\vspace{0.5cm}
   
    {\sc Nobuyuki Kenmochi}

 Faculty of Education, Chiba University, Chiba, Japan\\
  (E-mail: nobuyuki.kenmochi@gmail.com)\vspace{0.5cm}

   and\vspace{0.5cm}
     
  {\sc Marek Niezg\'odka}

 CNT Center, Cardinal Stefan Wyszynski University, Warsaw, Poland\\
  (E-mail: m.niezgodka@uksw.edu.pl)
 
\end{center}

\vspace{1cm}

\noindent
{\bf Abstract.} 
Variational inequalities, formulated on unknown--dependent convex 
sets, 
are called quasi-variational inequalities (QVI).
This paper is concerned with 
an abstract approach to a class of parabolic QVIs arising in many 
biochemical/mechanical problems. The approach is based on 
a compactness theorem for parabolic variational inequalities 
(cf. [10]). The prototype of our model for QVIs of
parabolic type is formulated in a reflexive Banach space 
as the sum of the time-derivative operator 
under unknown convex constraints and
a semimonotone operator, including a feedback system 
which selects a convex
constraint. The main objective of this work is 
to specify a class of unknown-state dependent convex constraints 
and to give
a precise formulation of QVIs.
\vspace{1cm}

\noindent
{\bf Keywords.} 
Variational inequalities, convex analysis, set--valued monotone operators, parabolic inequalities, superconductivity model.
\vspace{1cm}

\noindent
{\bf
AMS 2010 Subject Classification.} 34G25, 35G45, 35K51, 35K57, 35K59.

\vfill

\newpage
\noindent
\section{Introduction}
\label{sec:intro}

The theory of quasi-variational inequalities has been attractive since
its concept was formulated in [3]. Its development can be found in
[15, 21, 25, 29] for the stationary case and in [17, 23, 28] 
for the evolution case. 
But there are still many questions in the solvability of 
quasi-variational 
evolution inequalities, for instance, 
the establishment of compactness method 
from both of theoretical and numerical points of view. 
In this context, many mathematical models, arising in 
biochemical/mechanical problems 
(cf. [1, 2, 8, 11, 13, 14, 16, 19, 26, 27], 
have been discussed.

Our motivation for considering an abstract quasi--variational problem 
is as follows. Let $u_0$ be some initial state: 
e.g., an initial magnetic field and /or electric current, 
or the velocity of a flowing liquid, evolving in time to some 
unknown $u(t)$.
As a function of space, $u_0$ and $u(t)$  
belong to some Hilbert space $H$. The evolution is governed by  
some physical phenomenon that we describe by a system of equations,
like the Maxwell or the Navier--Stokes systems.
Here, we consider two general points.
\begin{enumerate}
\item Feedback. Within the described physical phenomenon, 
the evolution  may additionally depend on some parameter, 
that we call  $\theta$: e.g., the initial temperature, which 
is evolving according to its own physical law (e.g., the 
heat equation), but may be coupled to our unknown $u$ (the temperature
is strongly coupled to the magnetic state or to the concentration
of the substrate in the flow). 
\item Constraint. The very nature of the physical phenomena that our variable is undergoing may depend on space, on time, but also 
directly on the unknown itself.
In our examples, the magnetic material
passes to the superconductive state (zero or very low resistance)
when the magnetic field, the electric current and the temperature 
are low enough (actually, when some convex combination of these 
is low enough [5]); 
the flow velocity may depend on obstacles appearing due to
accumulation of some substrate that the flow is transporting ---
a phenomenon well known to cardiologists and their patients.
In the mathematical description, this means that the physical
equations governing the evolution are valid only in
some  subset of the space $H$, which is in general closed and convex.
This subset imposes a \emph{constraint} on the initial datum 
and on the unknown.   
Let us call this constraint set $\overline{K(\theta;t)}$.
\end{enumerate}
We are thus given a feedback operator for
$\theta$, which is taken from some metric space $\Theta$:
\begin{equation}
   \theta=\Lambda_{u_0} u,
\label{(1.4.b)}
\end{equation}
a family of convex sets $K(\theta;t)$ in $H$,
and we look for $u(t)\in {K(\theta;t)}$, which satisfies some
evolution equation. We formulate this last in an abstract way as
\begin{equation}
L_{u_0}(\theta;u)+A(u;u) \ni f. 
\label{(1.4.a)}
\end{equation}
Here, $L_{u_0}$ is a time derivative operator including the constraint
$u(t)\in {K(\theta;t)}$, and $A$ is a coercive semimonotone mapping,
possibly multivalued. 

Finally, the problem writes: 
given $u_0 \in H$, find $u(t)\in H$ such that
\begin{equation}
L_{u_0}(\theta;u)+A(u;u) \ni f,\qquad
   \theta=\Lambda_{u_0} u. 
\label{(1.4)}
\end{equation}
We give the functional framework for these abstract equations 
in Section~\ref{sec:frame}, the precise definitions 
of the convex sets  in Section~\ref{sec:Theta}, 
of the operators $L$ and $A$ in Sections~\ref{sec:L} and \ref{sec:pvi}, 
and we give two examples of these in Section~\ref{sec:app}. 
In particular, in Section~\ref{sec:superconductor} we give the concrete setting for the superconductivity type II problem, in a description inspired by
[1, 2, 19, 27]. As for the description of the flow problem, see our previous work [11, 12] and  references therein.

The main objective of this paper is to use the time-derivative 
operator $L_{u_0}(\theta; \cdot)$, by a 
systematic application of the semimonotone operator theory [22], 
to solve a class of parabolic quasi--variational inequalities of the form \eqref{(1.4)}. 

The paper is organized as follows. Section 2 introduces mathematical notation that we are going to use.  In Section \ref{sec:Theta}, we give 
the precise definition on the convex sets $K(\theta;t)$ and the assumptions on how they are related between them: we want some continuity of these with respect to $t$ and $\theta$. The continuity that we need for the sequel is given in Lemma~\ref{lem:2.1}, but we also give geometric conditions for this to occur: these are assumptions (A1)-(A3).
In Section \ref{sec:L}, we give the definition of the operator $L_{u_0}(\theta;\cdot)$, 
which is our "time derivative with constraint" operator.  
Fundamental properties of 
$L_{u_0}(\theta;\cdot)$ are discussed. In particular, geometric conditions (A1)-(A3) are shown to assure continuity of its graph (Theorem~\ref{th:3.2}).
Section~\ref{sec:pvi} is devoted to the operator $A$ and to a parabolic variational inclusion, intermediate for \eqref{(1.4.a)}: 
 $$ L_{u_0}(\theta; u)+A(v;u) \ni f. $$
We prove existence and uniqueness of its solution, 
as well as its continuous dependence upon the parameters 
$\theta$ and $v$. These are Theorems~\ref{th:4.1} and \ref{th:4.2}. 
In Section \ref{sec:pqvi}, the 
quasi-variational inequality \eqref{(1.4)} 
is formulated in detail and existence of solution is proved. 
Theorem~\ref{th:5.1}, stating existence of solution to~\eqref{(1.4)}, is the main theorem of this paper. We show immediately below that 
uniqueness cannot be expected. 
In Section \ref{sec:app} some applications, quasi-variational 
ordinary or partial differential inequalities, are given, showing that
our main assumptions  are not too much restrictive nor complicated to check.


\section{Functional framework and basic notation}
\label{sec:frame}

In general, for a (real) Banach space $X$, we denote by $X^*$ the 
dual space of $X$, by $|\cdot|_X$ and $|\cdot|_{X^*}$ their norms 
and by $\langle \cdot, \cdot \rangle_{X^*,X}$ the duality between 
$X^*$ and $X$. In the case when $X$ is a Hilbert space, the inner 
product may be denoted by $(\cdot,\cdot)_X$. 

For a proper, lower semicontinuous and convex function 
$\varphi$ on $X$, its effective domain $D(\varphi)$ and subdifferential
$\partial_{X,X^*}\varphi\: :\: X \to X^*$ 
are respectively defined by
 $$ D(\varphi):=\{z \in X~|~\varphi(z)<\infty\},$$
 $$\text{and }\quad \partial_{X,X^*}\varphi(z):=\{z^*\in X^*~|~\langle z^*, w-z\rangle
    _{X^*,X} \le
        \varphi(w) -\varphi(z), ~~\forall w \in X\},~\forall z \in X.$$
In the case when $X$ is a Hilbert space  
$\partial_{X,X^*} \varphi$ may be denoted simply by $\partial \varphi$.

Throughout this paper, let $H$ be a (real) Hilbert space 
and $V$ be a (real)
reflexive Banach space such that $V$ 
is dense and compactly embedded in $H$.
Then, identifying $H$ with the dual space $H^*$ of $H$, we have 
    $$V\subset H\subset V^*~~{\rm with~dense~and~compact~embeddings}. $$
We suppose that $V$ and $V^*$ are uniformly convex; hence the
duality mapping $F$ is singlevalued, continuous 
and strictly monotone from $V$ onto $V^*$.
Also, we write $\langle \cdot,\cdot \rangle$ for 
$\langle \cdot,\cdot \rangle_{V^*,V}.$ The positive numbers $p,~p'$
and $T$ are fixed so that
   $$ 2\le p <\infty,~~\frac 1p+\frac 1{p'}=1,~~0<T<\infty.$$ 
In this paper, we use the duality mapping $F: V\to V^*$ associated with the
gauge function $r \to r^{p-1},~r\geq 0$, namely
 $$ Fz \in V^*,~|Fz|_{V^*}=|z|_V^{p-1},~\langle Fz, z\rangle =|z|^p_V,~~
    \forall z \in V.$$

\section{Family of convex sets}
\label{sec:Theta}

\noindent
We are given here a complete metric space $\Theta$, with metric
$d_{\Theta}(\cdot,\cdot)$, consisting of  parameters $\theta$, 
and a family $\{K(\theta;t)\}_{\theta \in\Theta, t
\in [0,T]}$ of closed and convex subsets of~$V$ (see the previous section).

\begin{definition}
Wth the above assumptions, let
   $$ \psi^t(\theta;z)=\left\{
    \begin{array}{l}
 \displaystyle{I_{K(\theta;t)}(z)+\frac 1p|z|^p_V,~{\rm if}~ z \in V,}\\
 \displaystyle {\infty,   ~~~~~~~~~~~~~~~~~~~~{\rm if~} z \in H-V,}
    \end{array} \right.$$
for all $t \in [0,T]$ and for all $\theta\in \Theta$.
This is proper, l.s.c. and convex on $H$ and on $V$.
By the definition of subdifferential $\partial \psi^t(\theta;z)$ in $H$,
for any $z^* \in \partial \psi^t(\theta;z)$ we see that
  $$(z^*,w-z)_H \le \langle Fz, w-z\rangle,~~\forall w \in K(\theta;t).$$
\end{definition}

\subsection{Class of strong and weak parameters}

We aim at the precise formulation of the time-derivative 
$L_{u_0}(\theta;\cdot),~\theta \in \Theta$.
In this aim, we first  introduce a subset $\Theta_S$ of $\Theta$, 
which is the class of "strong parameters".
\begin{definition} $\Theta_S$
consists of all parameter $\theta \in \Theta$ such that 
$\{K(\theta;t)\}_{0\le t\le T}$ 
satisfies:
\begin{description}
\item{$(\Theta_S)$} for each positive number $r$ there are real-valued 
functions $a_{\theta,r} \in W^{1,2}(0,T)$ and $b_{\theta,r} \in W^{1,1}(0,T)$ 
having the following property that for any $s,~t \in [0,T]$ and any 
$z \in K(\theta;s)$ with $|z|_H \le r$ there is $\tilde z \in K(\theta;t)$ 
such that
 $$ |\tilde z-z|_H \le |a_{\theta,r}(t)-a_{\theta,r}(s)|(1+|z|_V^{\frac p2}),$$
and
 $$ |\tilde z|_V^p - |z|_V^p
         \le |b_{\theta,r}(t)-b_{\theta,r}(s)|(1+|z|_V^p).$$
\end{description}
\end{definition}

\begin{req}
It is well-known (cf. [18, 30]) that for a fixed 
$\theta\in \Theta_S$ condition $(\Theta_S)$
is sufficient for the Cauchy problem
 \begin{equation}
 u'(t) +\partial \psi^t(\theta;u(t)) \ni f(t) \text{~~in~}H, ~t\in [0,T],\qquad
     u(0)=u_0, \label{(2.1)}
\end{equation} 
to have a unique strong solution $u$ such that $u \in C([0,T];H)\cap 
L^p(0,T;V)$ with $\sqrt{t}u' \in L^2(0,T;H)$ and
$t \to |u(t)|^p_V$ is absolutely continuous on any compact interval of 
$(0,T]$ if 
$u_0 \in \overline{K(\theta;0)}$ and $f \in L^2(0,T;H)$; here 
$\partial \psi^t(\theta;\cdot)$ is the subdifferential of 
$\psi^t(\theta;\cdot)$
in $H$. In particular, if $u_0 \in K(\theta;0)$, then the solution $u$ 
of \eqref{(2.1)} belongs to $W^{1,2}(0,TH)$ and $t\to |u(t)|_V^p$ 
is absolutely continuous on 
$[0,T]$. This shows that  the set 
\begin{equation}
{\cal K}_0(\theta):=\{
   \eta \in L^p(0,T;V)~|~\eta' \in L^{p'}(0,T;V^*),
     ~\eta(t) \in K(\theta;t), \text{ a.e.~} t \in (0,T)\}; 
    \label{(1.2)}
\end{equation} 
which will be our class of test functions, is non-empty. 
\label{rem:2.1}
\end{req}

In the rest of this paper, we assume that the strong class  $\Theta_S$ is non-empty in $\Theta$.
\medskip

For a wider application of our theory 
we introduce the weak class $\Theta_W$:

\begin{definition}
\label{def:weak} $\Theta_W$, the weak class of parameters, is the
closure of $\Theta_S$ in $\Theta$.
\end{definition}
\subsection{More assumptions on the family of convex sets}

 We now specify the dependence of
$K(\theta;t)$ upon $\theta \in \Theta_W$ for each $t \in [0,T]$ by a family 
 $$ \{R_{\theta\bar \theta}(t) \in B(H)\cap B(V)~|~\theta,~\bar \theta \in 
\Theta_W,  ~t \in [0,T]\}$$
of bounded linear invertible operators in $H$ and $V$, where $B(H)$ 
(resp. $B(V)$) stands for the space of all bounded linear operators in $H$
 (resp. $V$), and by
a family $\{\sigma_{\theta\bar \theta,\varepsilon} \in C([0,T];V)~|~
\sigma'_{\theta,\bar\theta,\varepsilon} \in L^{p'}(0,T;V^*),~
\theta, \bar \theta \in \Theta_W,~ 0<  \varepsilon \le 1\}$ as follows.

\begin{description}
\item{(A1)} There is a positive constant $R_0$ such that
for all $\theta,~\bar \theta \in \Theta_W$
\begin{equation}
 |R_{\theta \bar \theta} -I|_{C([0,T];B(H))}
  +|R_{\theta \bar \theta} -I|_{C([0,T];B(V))} 
  +|R'_{\theta\bar \theta}|_{L^{p'}(0,T;B(H))}
    \le R_0d_\Theta(\theta,\bar \theta),
   \label{(2.2)}
\end{equation}
where $R'_{\theta\bar \theta}(t)= \frac d{dt}R_{\theta\bar \theta}(t)$ in 
$B(H)$. Moreover,
 for any $\theta,~\bar \theta \in \Theta_W$ we have
  \begin{equation}
  R^*_{\theta\bar \theta}(t)= R^{-1}_{\theta\bar \theta}(t)=
R_{\bar \theta \theta}(t),~~\forall
     \theta, \bar \theta \in \Theta_W,~\forall t \in [0,T], \label{(2.3)} 
\end{equation}   
namely the adjoint $R^*_{\theta\bar \theta}(t)$ of $R_{\theta\bar \theta}(t)$ 
coincides with
the inverse of $R_{\theta\bar \theta}(t)$ and it is $R_{\bar \theta \theta}(t)$ in $B(H)$.
\item{(A2)} There is a positive constant $\sigma_0$ such that
 $$ |\sigma_{\theta\bar \theta, \varepsilon}|_{C([0,T];V)}
    +|\sigma'_{\theta\bar \theta, \varepsilon}|_{L^{p'}(0,T;V^*)}  \le 
         \sigma_0(d_\Theta(\theta,\bar \theta)+\varepsilon),~~
   \forall \theta,~ \bar \theta \in \Theta_W,~\forall \varepsilon \in (0,1].$$
\item{(A3)} There is a continuous function $c_0(\cdot)$ on $[0,1]$ with
$c_0(0)=0$ satisfying the following property: 
 \begin{equation}
 \left \{
    \begin{array}{l}
  \forall \varepsilon \in (0,1],~\exists \delta_\varepsilon~
  \in (0,\varepsilon)~{\rm such~that~}\\[0.3cm]
  ~~~~\bar z={\cal F}_{\theta\bar \theta,\varepsilon}(t)z
    :=(1+c_0(\varepsilon))R_{\theta\bar \theta}(t)z +\sigma_{\theta\bar \theta,
    \varepsilon}(t)  \in K(\bar \theta;t),\\[0.2cm]
  ~~~~~~~~~~\forall z \in K(\theta;t),~\forall \theta,~\bar \theta \in 
   \Theta_W ~{\rm with~}
    d_\Theta(\theta,\bar \theta) \le \delta_\varepsilon, \forall t \in [0,T].
   \end{array} \right.
    \label{(2.4)}
    \end{equation}
\end{description}

\noindent
For $t \in [0,T],~\theta, \bar \theta \in \Theta_W$, (A3)  says 
that
$K(\theta;t)$ is mapped into $K(\bar \theta;t)$ by the operator 
$\bar z={\cal F}_{\theta\bar \theta, \varepsilon}(t)z$
given by \eqref{(2.4)}, 
which is a composition of rotation $R_{\theta\bar \theta}(t)$, 
contraction/expansion 
with modulus $(1+c_0(\varepsilon))$ close to 1 
and parallel transformation 
$\sigma_{\theta\bar \theta,\varepsilon}(t)$.

This type of transformation \eqref{(2.4)} was ealier introduced in [7] and [20].

\centerline{\includegraphics[height=0.36\textheight]{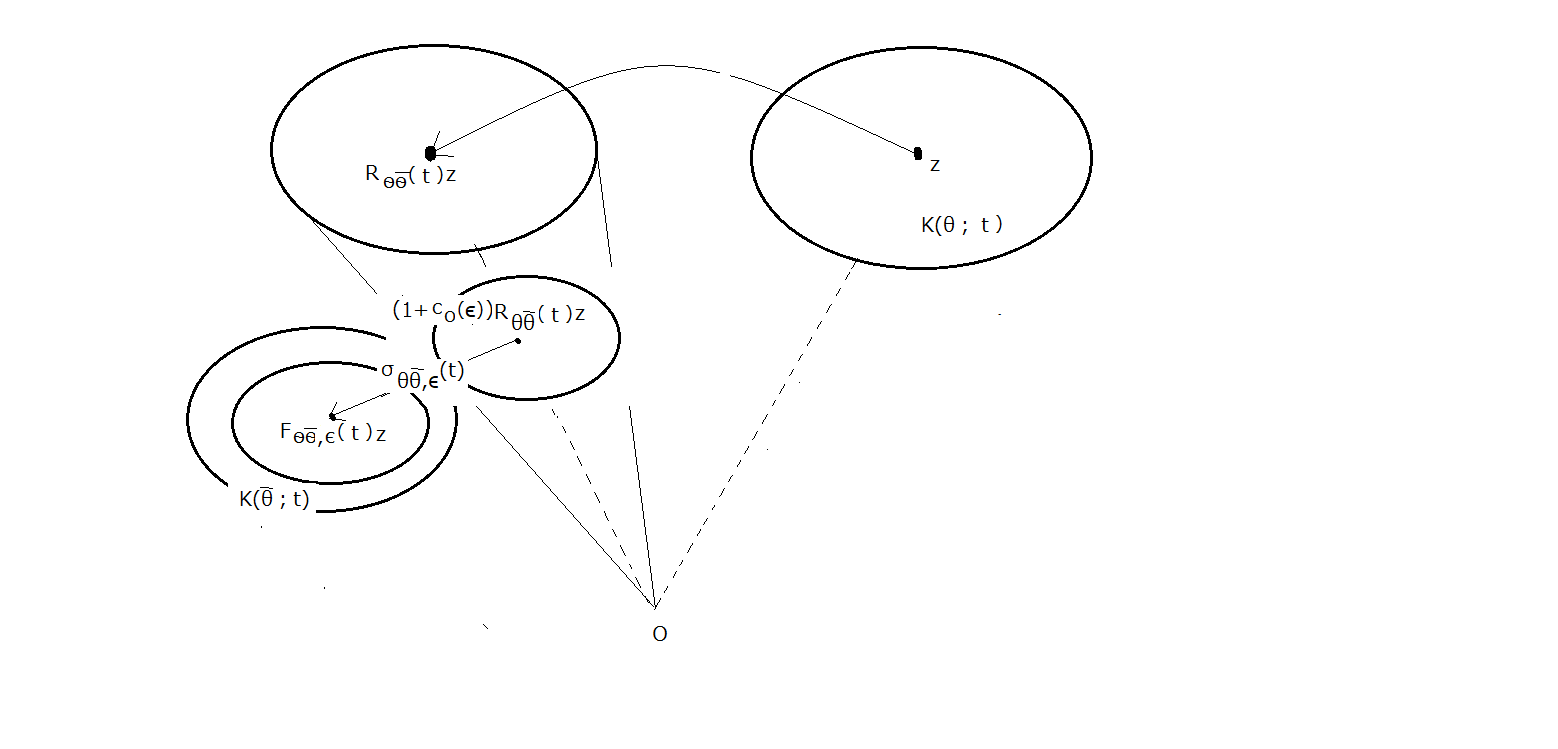}}


\subsection{Test functions for weak parameters}

As the class of test functions corresponding to parameter 
$\theta \in \Theta_W$, we employ
\begin{equation}
  {\cal K}(\theta):=\{\eta \in L^p(0,T;V)~|~\eta(t) \in K(\theta;t)~
  {\rm a.e.}~t \in (0,T)\} \label{(2.5)}
\end{equation}
The more narrow class ${\cal K}_0(\theta)$ of test functions corresponding to  parameter 
$\theta \in \Theta_S$ is defined  by \eqref{(1.2)}.
We see that ${\cal K}_0(\theta)$ is continuously embedded in $C([0,T];H)$.

The following lemma states that for the weak parameters, test functions from the above set subject to our transformation 
${\cal F}_{\theta\bar \theta, \varepsilon}(t)$, converge to the original function when $\varepsilon$ tends to zero and $\bar\theta$ tends to $\theta$. This is a condition for a kind of Mosco convergence  [24].
\label{lem:2.1}
\begin{lem}
 Let $\theta,~\bar \theta \in \Theta_W$ and $\varepsilon$ 
be any small positive number. For each $\eta \in {\cal K}(\theta)$, put
$\eta_{\theta\bar \theta,\varepsilon}(t):={\cal F}_{\theta\bar \theta,
\varepsilon}(t)\eta(t)$ for a.e. $t \in (0,T)$. Then, 
$\eta_{\theta\bar \theta,\varepsilon} \in {\cal K}(\bar \theta)$ and 
$ \eta_{\theta\bar \theta,\varepsilon} \to \eta~{\it in~}L^p(0,T;V)$
as $\bar \theta \to \theta$ in $\Theta$ and $\varepsilon \downarrow 0$,
and moreover,
if $\eta \in {\cal K}_0(\theta)$, then
\begin{equation}
\eta_{\theta\bar \theta,\varepsilon} \to \eta~{\it in~}C([0,T];H),
   ~~\eta'_{\theta\bar \theta,\varepsilon} \to \eta'
    ~{\it in~}L^{p'}(0,T;V^*),\label{(2.8)}
    \end{equation}
as $\bar \theta \to \theta$ in $\Theta$ and $\varepsilon \downarrow 0.$
 \end{lem}
\noindent
{\bf Proof.} Let $\eta \in {\cal K}(\theta)$. We see from \eqref{(2.2)}-\eqref{(2.3)} of (A1) 
and (A2) that 
$\eta_{\theta\bar \theta,\varepsilon}(t) \in K(\bar \theta; t)$ for a.e.  
$t \in (0,T)$ and
\begin{equation}
\begin{array}{l}
   \displaystyle{
  |\eta_{\theta\bar \theta,\varepsilon}(t)-\eta(t)|_V =|(R_{\theta\bar \theta}
   (t)-I)\eta(t)+c_0(\varepsilon)R_{\theta\bar \theta}(t)\eta(t) 
   +\sigma_{\theta\bar \theta,\varepsilon}(t)|_V}\\[0.3cm]
\displaystyle {~~~~~~~~~~~~~~~~~~~~~\le (R_0d_\Theta(\theta,\bar \theta)
    +|c_0(\varepsilon)|
   |R_{\theta\bar\theta}|_{C([0,T];B(V))})|\eta(t)|_V}\\[0.3cm]
\displaystyle{~~~~~~~~~~~~~~~~~~~~~~~~~
  +\sigma_0(d_\Theta(\theta,\bar \theta)+\varepsilon)}
    \end{array} \label{(2.9)} 
\end{equation} 
whence $\eta_{\theta\bar \theta,\varepsilon} \to \eta$ in $L^p(0,T;V)$ as 
 $\bar \theta \to \theta$ in $\Theta$ and $\varepsilon \downarrow 0$.

Next, let $\eta\in {\cal K}_0(\theta)$. Then $\eta_{\theta\bar\theta,
\varepsilon},~\eta \in C([0,T];H)$ and
in a similar way to \eqref{(2.9)}
  $$ 
  |\eta_{\theta\bar \theta,\varepsilon}(t)-\eta(t)|_H 
   \le (R_0d_\Theta(\theta,\bar \theta)
    +|c_0(\varepsilon)|
   |R_{\theta\bar\theta}|_{C([0,T];B(H))})|\eta(t)|_H
  +\sigma_0(d_\Theta(\theta,\bar \theta)+\varepsilon)$$
which implies that $\eta_{\theta\bar\theta,\varepsilon}$ converges to $\eta$
in $H$ uniformly on $[0,T]$, namely in $C([0,T];H)$ as $\bar\theta \to \theta$
in $\Theta$ and $\varepsilon \to 0$. 
As to the time derivative of 
$\eta_{\theta\bar \theta, \varepsilon}$, 
for every $z \in V$ we observe that
  \begin{equation}
   \begin{array}{l}
  \displaystyle{ \frac d{dt}\langle R_{\theta\bar \theta}(t)\eta(t), z\rangle
  =\frac d{dt}(R_{\theta\bar \theta}(t)\eta(t), z)_H= \frac d{dt}(\eta(t), 
    R_{\bar \theta \theta}(t)z)_H}\\[0.3cm]
  \displaystyle{~~~~~~~~~~~~~~~~~~~~~~
    = \langle \eta'(t), R_{\bar \theta \theta}(t)z \rangle +(\eta(t),
   R'_{\bar \theta \theta}(t)z)_H }
  \end{array} \label{(2.10)}
\end{equation}  
for a.e. $t \in (0,T)$. This shows that $R_{\theta\bar \theta}(t)\eta(t)$ is 
(strongly) differentiable in $t$ a.e. on $(0,T)$, since the last term of \eqref{(2.10)}
is the uniform limit of 
 $$\lim_{\Delta t\to 0} \frac 1{\Delta t}
    \left\{\langle \eta(t+\Delta t)-\eta (t),
  R_{\bar \theta\theta}(t)z\rangle 
  + (\eta(t), (R_{\bar\theta\theta}(t+\Delta t)-R_{\bar\theta\theta}(t))z)_H
    \right\}$$
with respect to $z$ with $|z|_V\le 1$ at a.e. 
$t \in (0,T)$.
Also, it follows from \eqref{(2.10)} that if $|z|_V \le 1$, then 
  \begin{eqnarray*}
    \left|\frac d{dt}R_{\theta\bar\theta}(t)\eta(t)\right|_{V^*}
    &\le& \sup_{|z|_V\le 1}\left\{ 
     |\eta'(t)|_{V^*}|R_{\bar\theta\theta}(t)|_{B(V)}|z|_V + |\eta(t)|_H
     |R'_{\bar\theta\theta}(t)|_{B(H)}\cdot c_V|z|_V\right\}\\[0.2cm]
    &\le& |\eta'(t)|_{V^*}|R_{\bar\theta\theta}(t)|_{B(V)} + c_V|\eta(t)|_H
     |R'_{\bar\theta\theta}(t)|_{B(H)}\\
    &\le& |\eta'(t)|_{V^*}|R_{\bar\theta\theta}|_{C([0,T];B(V))} 
    + c_V|\eta|_{C([0,T:H)}|R'_{\bar\theta\theta}(t)|_{B(H)},
  \end{eqnarray*}
where $c_V$ is a positive constant satisfying that $|w|_H \le c_V|w|_V$ for
all $w \in V$. This shows that  $\frac d{dt}R_{\theta\bar \theta}\eta \in 
L^{p'}(0,T;V^*)$ and hence $\eta'_{\theta\bar \theta,\varepsilon}=
(1+c_0(\varepsilon))
\frac d{dt}R_{\theta\bar \theta}\eta+\sigma'_{\theta\bar\theta,\varepsilon} 
\in L^{p'}(0,T;V^*)$.
Besides, we see from \eqref{(2.10)} that
  $$ \langle \frac d{dt}R_{\theta\bar \theta}(t)\eta(t) - \eta'(t), z \rangle
     =\langle \eta'(t), (R_{\bar\theta\theta}(t)-I)z\rangle
       +(\eta(t),R'_{\bar \theta \theta}(t)z)_H, $$
which shows that
\begin{equation}
\left|\frac d{dt}R_{\theta\bar \theta}(t)\eta(t) - \eta'(t) \right|_{V^*} 
    \le |\eta'(t)|_{V^*}|R_{\bar\theta\theta}-I|_{C([0,T];B(V))}
  +c_V|\eta|_{C([0,T];H)}
    |R'_{\bar\theta\theta}(t)|_{B(H)}. \label{(2.11)}
\end{equation} 
Since 
 $$\eta'_{\theta\bar \theta,\varepsilon}(t)-\eta'(t)
   =\frac d{dt}R_{\theta\bar \theta}(t)\eta(t) -\eta'(t) +c_0(\varepsilon)
    \frac d{dt}R_{\theta\bar \theta}(t)\eta(t) 
    +\sigma'_{\theta\bar \theta,\varepsilon}(t),$$
it results from \eqref{(2.11)} that $\eta'_{\theta\bar \theta,\varepsilon}\to \eta'$
in $L^{p'}(0,T;V^*)$ as $\bar \theta \to \theta$ in $\Theta$ and 
$\varepsilon \downarrow 0$. Thus \eqref{(2.8)} is obtained. \hfill $\Box$ 

\section{Time-derivative operator $L_{u_0}(\theta; \cdot)$}
\label{sec:L}

\begin{definition}
Let $\theta \in \Theta_W$ and 
$u_0 \in \overline{K(\theta;0)}$. We define the operator
$L_{u_0}(\theta;\cdot)$ by: $g \in L_{u_0}(\theta;u)$ 
if and only if 
\begin{equation}
\begin{array}{l}
  \displaystyle{u \in {\cal K}(\theta),
   ~g\in L^{p'}(0,T;V^*),}\\[0.3cm]
  \displaystyle{\int_0^T \langle \eta'-g,u-\eta \rangle dt 
   \le \frac 12|u_0-\eta(0)|^2_H,~~\forall \eta \in {\cal K}_0(\theta),}
    \end{array}  \label{(2.7)}
\end{equation}
where ${\cal K}(\theta)$ and ${\cal K}_0(\theta)$ are the class of test functions given by \eqref{(2.5)} and \eqref{(1.2)}, respectively.
\end{definition}

The fundamental properties of $L_{u_0}(\theta;\cdot)$ as well as its 
characterization are derived from the proposition stated below.

\begin{prop}
Suppose that
$(\Theta_S)$ and
(A1)-(A3) are fulfilled. Let $\varepsilon$ be any small positive number and
 $\theta_i \in \Theta_S,~i=1,2,$ such that $d_\Theta(\theta_1,\theta_2)
 < \varepsilon$.
Also, let $M$ be any positive number, 
$u_{0i} \in \overline {K(\theta_i;0)}$, 
and $f_i \in L^2(0,T;H), ~i=1,2,$ such that
\begin{equation}
   \sum_{i=1}^2 \left\{|f_i|_{L^{p'}(0,T;V^*)}+|u_{0i}|_H\right\}
     \le M. \label{(2.12)}
\end{equation} 
Then, for the strong solutions $u_i$  of
\begin{equation}
u'_i(t)+\partial \psi^t(\theta_i; u_i(t)) \ni f_i(t)~~{\it in~}H,
   ~{\it a.e.~}t \in (0,T),~~u_i(0)=u_{0i},~i=1,2,
    \label{(2.13)}
\end{equation} 
we have that for all $s,~t \in [0,T]$ with $s\le t$ 
\begin{equation}
\begin{split}
{\frac 12 |u_1(t)-u_2(t)|^2_H 
     +\int_0^t \langle Fu_1(\tau)-Fu_2(\tau), 
     u_1(\tau)-u_2(\tau) \rangle d\tau} \\[0.3cm]
   \displaystyle{\le \frac 12 |u_{1}(s)-u_{2}(s)|^2_H 
   +\int_s^t (f_1(\tau)-f_2(\tau), 
   u_1(\tau)-u_2(\tau) )_H d\tau + C^*_0(M) C^*_1(\varepsilon)},
   \end{split}
 \label{(2.14)}  
\end{equation}
where $C^*_0(M)$ is a positive constant depending only on $M>0$ and 
$C^*_1(\varepsilon)$ is a positive continuous 
function of $\varepsilon$ satisfying $C^*_1(\varepsilon) \to 0$ as 
$\varepsilon \downarrow 0$.
\label{prop:2.1}
\end{prop}

Under the same assumptions and same notation as in Proposition \ref{prop:2.1},
as was noted in Remark \ref{rem:2.1}, problem \eqref{(2.1)} admits a unique solution
$u_i \in C([0,T];H)$ with $\sqrt{t}u'_i \in L^2(0,T;H)$. Put $u^*_i(\tau) :=
f_i(\tau)-u'_i(\tau)\in \partial \psi^{\tau}(\theta_i;u_i(\tau))$ for a.e. 
$\tau \in (0,T)$. Since $\psi^{\tau}(\theta_i;\cdot)$ is coercive, namely
$\psi^{\tau}(\theta_i;z) \geq \frac 1p|z|^p_V$ for all $z \in V$, we note 
from the usual energy estimate for problem \eqref{(2.1)} that 
  \begin{equation}
  |u_i|_{C([0,T];H)}+|u_i|_{L^p(0,T;V)} \le M_i:=
   M_i(|f_i|_{L^{p'}(0,T;V^*)},~|u_{0i}|_H), ~i=1,2, 
   \label{(2.15)}
\end{equation}
where 
$M_i(\cdot,\cdot)$ is a non-negative and non-decreasing function on
${\rm R}^2$.

We take the inner product between the both sides of 
$u'_1(\tau) +u^*_1(\tau)=f_1(\tau)$ and 
$u_1(\tau)-{\cal F}_{\theta_2\theta_1,\varepsilon}(\tau)u_2(\tau)$ to obtain 
by the definition of the subdifferential of 
$\partial \psi^{\tau}(\theta_1;\cdot)$
  $$ (u'_1(\tau),u_1(\tau)-{\cal F}_{\theta_2\theta_1,\varepsilon}(\tau)
   u_2(\tau))_H
     +\langle Fu_1(\tau),u_1(\tau)-{\cal F}_{\theta_2 \theta_1,\varepsilon}
   (\tau)u_2(\tau)\rangle$$
$$ \le (f_1(\tau),u_1(\tau)-{\cal F}_{\theta_2\theta_1,\varepsilon}(\tau)
  u_2(\tau))_H,~~{\rm a.e.~}\tau \in (0,T).\hspace{3.5cm} $$
Here we observe that
\begin{equation}
\begin{array}{l}
  \displaystyle{(u'_1(\tau),u_1(\tau)-u_2(\tau))_H
      +\langle Fu_1(\tau),u_1(\tau)-u_2(\tau)\rangle }\\[0.3cm]
   \displaystyle{ \le (f_1(\tau),u_1(\tau)-u_2(\tau))_H 
     + \Gamma_{\theta_1\theta_2,\varepsilon}(\tau)
     +\tilde \Gamma_{\theta_1\theta_2,\varepsilon}(\tau)}
  \end{array} 
      \label{(2.16)}
\end{equation} 
with
 \begin{equation}
  \begin{array}{l}
  \displaystyle{\Gamma_{\theta_1\theta_2,\varepsilon}(\tau)
    := 
    (u'_1(\tau), (R_{\theta_2\theta_1}(\tau)-I)u_2(\tau))_H +c_0(\varepsilon)
   (u'_1(\tau), R_{\theta_2\theta_1}(\tau)u_2(\tau))_H }\\[0.3cm]
    \displaystyle{~~~~~~~~+\frac d{d\tau}
    (u_1(\tau),\sigma_{\theta_2\theta_1,\varepsilon}(\tau))_H}\\[0.3cm]
  \displaystyle{ \tilde\Gamma_{\theta_1\theta_2,\varepsilon}(\tau)
    := (f_1(\tau), (I-R_{\theta_2\theta_1}(\tau))u_2(\tau)-c_0(\varepsilon)
 R_{\theta_2\theta_1}(\tau)u_2(\tau)-\sigma_{\theta_2\theta_1,\varepsilon}
  (\tau))_H}
  \\[0.3cm]
 \displaystyle{~~~~~~~~-\langle \sigma'_{\theta_2\theta_1,
  \varepsilon}(\tau), u_1(\tau)\rangle - \langle Fu_1(\tau), u_2(\tau)- 
   {\cal F}_{\theta_2\theta_1,\varepsilon}(\tau)u_2(\tau)\rangle.}
       \end{array} 
   \label{(2.17)}
\end{equation} 
Similarly, by exchanging parameters $\theta_1$ and $\theta_2$ we have
\begin{equation} \begin{array}{l}
 \displaystyle{ (u'_2(\tau),u_2(\tau)-u_1(\tau))_H
      +\langle Fu_2(\tau),u_2(\tau)-u_1(\tau)\rangle }\\[0.3cm]
    \displaystyle{\le (f_2(\tau),u_2(\tau)-u_1(\tau))_H 
   + \Gamma_{\theta_2\theta_1,\varepsilon}(\tau)
   + \tilde \Gamma_{\theta_2\theta_1,\varepsilon}(\tau),}
      \end{array} \label{(2.18)}
\end{equation}
Adding \eqref{(2.16)} and \eqref{(2.18)} yields that
  $$ \begin{array}{l}
   \displaystyle{~~~\frac 12|u_1(t)-u_2(t)|^2_H 
       + \int_0^t \langle Fu_1-Fu_2, u_1-u_2\rangle d\tau}\\[0.3cm] 
   \displaystyle{ \le \frac 12 |u_{01}-u_{02}|^2_H 
      + \int_0^t(f_1-f_2,u_1-u_2)_H d\tau} \\[0.2cm]
   \displaystyle{ ~~~~~+\int_0^t(\Gamma_{\theta_2\theta_1,\varepsilon}
    +\Gamma_{\theta_1\theta_2,\varepsilon}
     +\tilde \Gamma_{\theta_2\theta_1,\varepsilon}+\tilde 
     \Gamma_{\theta_1\theta_2,\varepsilon})d \tau,
      ~~\forall t \in [0,T].}
    \end{array}  $$
Now, recall the rearrangement of $\Gamma_{\theta_2\theta_1,\varepsilon}, 
~\tilde \Gamma_{\theta_2\theta_1,\varepsilon}$ and 
$\Gamma_{\theta_1\theta_2,\varepsilon},
~\tilde \Gamma_{\theta_1\theta_2,\varepsilon}$, which can be derived from 
the assumptions
(A1), (A2) and (A3) on $R_{\theta\bar \theta}(t)$ and 
$\sigma_{\theta\bar \theta,\varepsilon}$
together with \eqref{(2.17)}. 

\begin{lem}
\label{lem:2.2} 
([20; Lemma 3.2]) We have: for a.e. $\tau \in (0,T)$,
 \begin{eqnarray*}
   && \Gamma_{\theta_1\theta_2,\varepsilon}(\tau)
   +\Gamma_{\theta_2\theta_1, \varepsilon}(\tau) \\
   &=&\frac d{d\tau}(u_1(\tau),
   (R_{\theta_2\theta_1}(\tau)-I)u_2(\tau))_H
  +c_0(\varepsilon)\frac d{d\tau} (u_1(\tau),R_{\theta_2\theta_1}    
    (\tau)u_2(\tau))_H \\
   && +\frac d{d\tau}\{(u_1(\tau),\sigma_{\theta_2\theta_1,
     \varepsilon}(\tau))_H+(u_2(\tau), \sigma_{\theta_1\theta_2,\varepsilon}
     (\tau))_H \} \\
   && -(1+c_0(\varepsilon))(u_1(\tau), 
       R'_{\theta_2\theta_1}(\tau)u_2(\tau))_H, 
  \end{eqnarray*}  
 \begin{eqnarray*}
   && |\tilde\Gamma_{\theta_1\theta_2,\varepsilon}(\tau)|
   +|\tilde\Gamma_{\theta_2\theta_1, \varepsilon}(\tau)| \\
   &\le& |f_1(\tau)|_{V^*}\{|I-R_{\theta_2\theta_1}
    (\tau)|_{B(V)}
    |u_2(\tau)|_V+|c_0(\varepsilon)||u_2(\tau)|_V +|\sigma_{\theta_2\theta_1,
     \varepsilon}(\tau)|_V\} \\
   & & +|u_1(\tau)|_V|\sigma'_{\theta_2\theta_1,
    \varepsilon}(\tau)|_{V^*} + |u_1(\tau)|^{p-1}_V|{\cal F}_{\theta_2\theta_1,
    \varepsilon}(\tau)u_2(\tau)-u_2(\tau)|_V \\
   && +|f_2(\tau)|_{V^*}\{|I-R_{\theta_1\theta_2}
   (\tau)|_{B(V)}|u_1(\tau)|_V+|c_0(\varepsilon)||u_1(\tau)|_V 
    +|\sigma_{\theta_1\theta_2,
     \varepsilon}(\tau)|_V\} \\
   && +|u_2(\tau)|_V|\sigma'_{\theta_1\theta_2,
    \varepsilon}(\tau)|_{V^*} +|u_2(\tau)|^{p-1}_V|{\cal F}_{\theta_1\theta_2,
    \varepsilon}(\tau)u_1(\tau)-u_1(\tau)|_V
     \end{eqnarray*} 
and 
  \begin{eqnarray*}  
&& |u_1(\tau)|^{p-1}_V |{\cal F}_{\theta_2\theta_1,\varepsilon}u_2(\tau)
     -u_2(\tau)|_V + |u_2(\tau)|^{p-1}_V|{\cal F}_{\theta_1\theta_2,
    \varepsilon}u_1(\tau)-u_1(\tau)\rangle|_V\\
&\le&|u_1(\tau)|^{p-1}_V(|I-R_{\theta_2\theta_1}(\tau)|_{B(V)}|u_2(\tau)|_V+
   |c_0(\varepsilon)||R_{\theta_2\theta_1}(\tau)|_{B(V)}|u_2(\tau)|_V
   +|\sigma_{\theta_2\theta_1,\varepsilon}(\tau)|_V)\\
&&\hspace{-0.4cm}+|u_2(\tau)|^{p-1}_V(|I-R_{\theta_1\theta_2}(\tau)|_{B(V)}|u_1(\tau))|_V+
   |c_0(\varepsilon)||R_{\theta_1\theta_2}(\tau)|_{B(V)}|u_1(\tau)|_V
   +|\sigma_{\theta_1\theta_2,\varepsilon}(\tau)|_V).
  \end{eqnarray*} 
\end{lem}
\noindent
\hspace*{-0.15cm}{\bf Proof of Proposition \ref{prop:2.1}.} 
Let $M$ and $M_i$ be the same constants as 
in \eqref{(2.12)} and \eqref{(2.15)}, and let $\varepsilon$ be any small positive number. Then,
by using (A1)-(A3), we have 
on account of Lemma \ref{lem:2.2} that for all
$s, ~t \in [0,T],~s\le t$.
  \begin{eqnarray*}
 && \left| \int_s^t (\Gamma_{\theta_2\theta_1,\varepsilon}
  +\Gamma_{\theta_1\theta_2,\varepsilon}) d\tau\right|\\
 & \le& 2M_1M_2|I-R_{\theta_2\theta_1}|_{C([0,T];B(H))}
 +2|c_0(\varepsilon)|M_1M_2
 +2M_1 |\sigma_{\theta_2\theta_1,\varepsilon}|_{C([0,T];H)}\\
 && +2M_2|\sigma_{\theta_1\theta_2,\varepsilon}|_{C([0,T];H)}
  + (1+|c_0(\varepsilon)|)M_1M_2 T^{\frac 1p}
    |R'_{\theta_2\theta_1}|_{L^{p'}(0,T;B(H))}\\
 &\le& 2d_\Theta(\theta_1,\theta_2)+2|c_0(\varepsilon)| M_1M_2
   +2(M_1+M_2)c_V \sigma_0(d_\Theta(\theta_1,\theta_2)+\varepsilon)\\
  && +2M_1M_2T^{\frac 1p}R_0 d_\Theta(\theta_1,\theta_2);
  \end{eqnarray*}
we used above $|R_{\theta_i\theta_j}(\tau)|_{B(H)}=1$ for all $\tau$ and
the inequality $|z|_H \le c_V|z|_V$ for all $z \in V$ with an embedding 
constant $c_V>0$ to have
$|\sigma_{\theta_2\theta_1,\varepsilon}|_{C([0,T];H)} \le c_V
   |\sigma_{\theta_2\theta_1,\varepsilon}|_{C([0,T];V)}$.
Similarly,
  \begin{eqnarray*}
   && \int_0^T \{|\tilde\Gamma_{\theta_2\theta_1,\varepsilon}|
  +|\tilde\Gamma_{\theta_1\theta_2,\varepsilon}| \} d\tau \\
  &\le& M(M_1+M_2)(|I-R_{\theta_2\theta_1}|_{C([0,T];B(V))}
     +|I-R_{\theta_2\theta_1}|_{C([0,T];B(V))}) \\
  && +|c_0(\varepsilon)|T^{\frac 1{p'}}(M_1+M_2) 
  + M(|\sigma_{\theta_2\theta_1,\varepsilon}|_{C([0,T];V)}
   + |\sigma_{\theta_1\theta_2,\varepsilon}|_{C([0,T];V)}) \\
  && +M_1|\sigma'_{\theta_2\theta_1,\varepsilon}|_{L^{p'}(0,T;V^*)} 
   +M_2|\sigma'_{\theta_1\theta_2,\varepsilon}|_{L^{p'}(0,T;V^*)} \\
  && +(M_1^{\frac p{p'}}+M_2^{\frac p{p'}}) 
   (|I-R_{\theta_2\theta_1}|_{C([0,T];B(V))}
   +|I-R_{\theta_1\theta_2}|_{C([0,T];B(V))}+2|c_0(\varepsilon)|)\\
    && + M_1^{\frac p{p'}}|\sigma_{\theta_1\theta_2,\varepsilon}|_{L^p(0,T;V)}
    + M_2^{\frac p{p'}}|\sigma_{\theta_2\theta_1,\varepsilon}|_{L^p(0,T;V)}
   \end{eqnarray*}
    \begin{eqnarray*}
 & \le& 2M(M_1+M_2)R_0d_\Theta(\theta_1,\theta_2)+ T^{\frac 1{p'}}(M_1+M_2)
   |c_0(\varepsilon)|+ 2M\sigma_0(d_\Theta(\theta_1\theta_2)+\varepsilon)\\
 && +(M_1+M_2)\sigma_0(d_\Theta(\theta_1,\theta_2)+\varepsilon)
   + (M_1^{\frac p{p'}}+M_2^{\frac p{p'}}(2R_0d_\Theta(\theta_1,\theta_2)+
    |c_0(\varepsilon)|)\\
 &&   + (M_1^{\frac p{p'}}+M_2^{\frac p{p'}} )T^{\frac 1p}\sigma_0
    (d_\Theta(\theta_1,\theta_2)+\varepsilon). 
  \end{eqnarray*}
By these estimates we can find some constant
$C^*_0(M)$ and function $C^*_1(\varepsilon)$ of $\varepsilon \in
(0,1]$ so that \eqref{(2.14)}
holds, for instance,
$ C^*_0(M) :=(M_1+M_2)(4c_V\sigma_0+2MR_0+T^{\frac 1{p'}} +2\sigma_0
     +T^{\frac 1{p'}}) + 2R_0(M_1^p+M_2^p)+ 
4(M_1^{\frac p{p'}}+M_2^{\frac p{p'}})
(2R_0+T^{\frac 1p}\sigma_0)+2M_1M_2(1+T^{\frac 1p}R_0) +2$ and 
 $C^*_1(\varepsilon):=\varepsilon+ |c_0(\varepsilon)|$.
\hspace*{1cm}\hfill $\Box$

\label{sec:prop}

Now, we investigate some fundamental properties of 
$L_{u_0}(\theta;\cdot)$ which are derived from Proposition \ref{prop:2.1}.

\begin{lem}
\label{lem:3.1}
Let $\theta \in \Theta_W$, $u_0 \in 
\overline{K(\theta;0)}$ and $f \in L^{p'}(0,T;V^*)$, and let $\{\theta_n\}
\subset \Theta_S$, $\{u_{0n}\} \subset H$ and $\{f_n\} \subset L^2(0,T;H)$
such that $u_{0n} \in \overline{K(\theta_n;0)}$,
\begin{equation}
\theta_n \to \theta~{\it in~}\Theta,~u_{0n} \to u_0~{\it in~} H,~
     f_n \to f~{\it in~} L^{p'}(0,T;V^*)~
    ~({\it as~} n \to \infty). \label{(3.1)}
\end{equation} 
Then the strong solution $u_n$ of 
   \begin{equation}
   u'_n(t)+\partial \psi^t(\theta_n;u_n(t))\ni f_n(t)~{\it in~}H,~ a.e.~
     t\in (0,T),~u_n(0)=u_{0n}, \label{(3.2)}
\end{equation}    
converges in $C([0,T];H)\cap L^p(0,T;V)$ to a solution $u$ of
  \begin{equation}
  L_{u_0}(\theta; u) +Fu \ni f~~{\it in~}L^{p'}(0,T; V^*). \label{(3.3)}
\end{equation}   
\end{lem}
\noindent
{\bf Proof.} 
Given any small $\varepsilon>0$ and a constant $M>0$ satisfying
$ |f|_{L^{p'}(0,T;V^*)} +|u_0|_H +1 \le \frac M2 $,
choose a positive integer $N_\varepsilon$
such that $d_\Theta(\theta_n, \theta) \le \frac {\varepsilon}2$ and
$|f_n|_{L^{p'}(0,T;V^*)}+|u_{0n}|_H \le \frac M2$ for all 
$n \geq N_\varepsilon$. Then it follows from Proposition \ref{prop:2.1} with \eqref{(3.1)} that
\begin{equation}
\begin{array}{l}
   \displaystyle{~~~\frac 12 |u_n(t)-u_m(t)|^2_H 
     +\int_0^t \langle Fu_n-Fu_m, 
     u_n-u_m  \rangle d\tau} \\[0.3cm]
   \displaystyle{\le \frac 12 |u_{0n}-u_{0m}|^2_H 
   +\int_0^t (f_n-f_m, 
   u_n-u_m )_H d\tau + C^*_0(M) C^*_1(\varepsilon),}\\[0.3cm]~~
   \displaystyle{\hspace{5cm}  \forall t \in [0,T],
     ~~\forall n,~m \geq N_\varepsilon.}
    \end{array} \label{(3.4)}
\end{equation} 
which implies that $u_n \to u$ in $C([0,T];H)$ as $n \to \infty$ and 
$\varepsilon \to 0$ for a certain
function $u\in C([0,T];H) \cap L^p(0,T;V)$ and
  \begin{equation}
   \lim_{n,m\to \infty} \int_0^T \langle Fu_n-Fu_m, u_n-u_m  \rangle d\tau 
      =0. \label{(3.5)}
\end{equation}  
By the uniform convexity of $V$, it follows from \eqref{(3.5)} that 
$u_n \to u$ in $L^p(0,T;V)$. Since for large $n$, small $\varepsilon>0$ and 
any $\eta \in K_0(\theta)$, with notation $\eta_{n,\varepsilon}(t):=
{\cal F}_{\theta \theta_n,\varepsilon}(t)\eta(t)$
we see from \eqref{(3.2)} that for any interval $[s,t] \subset [0,T]$
  $$ \int_s^t (u'_n, u_n- \eta_{n,\varepsilon})_H d\tau
  +\int_s^t \langle Fu_n, u_n - \eta_{n,\varepsilon}\rangle d\tau
  \le \int_s^t \langle f_n, u_n -\eta_{n,\varepsilon}\rangle d\tau
  $$
from which we obtain by integration by parts
   \begin{eqnarray*}   
  & &\int_s^t \langle \eta'_{n,\varepsilon}, u_n- \eta_{n,\varepsilon}\rangle
    d\tau +\frac 12|u_n(t)-\eta_{n,\varepsilon}(t)|^2_H
     +\int_0^t \langle Fu_n, u_n - \eta_{n,\varepsilon}\rangle d\tau \\
  &\le& \frac 12 |u_n(s)-\eta_{n,\varepsilon}(s)|^2_H +\int_s^t
     \langle f_n, u_n -    \eta_{n,\varepsilon}\rangle d\tau.
   \end{eqnarray*}
Now, passing to the limit in this inequality as $\theta_n \to \theta$ 
in $\Theta$
and $\varepsilon \downarrow 0$ by using Lemma \ref{lem:2.1}, we conclude that 
 \begin{equation}
  \int_s^t \langle \eta' -f+Fu, u- \eta\rangle d\tau 
    +\frac 12|u(t)-\eta(t)|^2_H
  \le \frac 12 |u(s)-\eta(s)|^2_H, ~~\forall \eta \in {\cal K}_0(\theta).
\label{(3.6)} 
\end{equation} 
This implies $ f - Fu \in L_{u_0}(\theta;u)$. Thus \eqref{(3.3)} is obtained.
\hfill $\Box$ 

\begin{cor}
\label{cor:3.1} 
Let $\theta_i \in \Theta_W$, $u_{0i} \in \overline
{K(\theta_i;0)}$ and $f_i \in L^{p'}(0,T;V^*)$, $i=1,2$, and let $u_i$ be a
solution of $L_{u_{0i}}(\theta_i;u_i)+Fu_i \ni f_i$ in $L^{p'}(0,T;V^*)$, 
$i=1,2$. Let $M>0$ and $\varepsilon$ be any positive numbers such that 
  $$\sum_{i=1}^2 \left\{|f_i|_{L^{p'}(0,T;V^*)}+|u_{0i}|_H\right\} \le M,~~ 
             d_\Theta(\theta_1,\theta_2) \le \varepsilon.$$
Then, for any $s,~t \in [0,T],~s \le t$,
\begin{equation}
\begin{array}{l}
 \displaystyle{\frac 12|u_1(t)-u_2(t)|^2_H +\int_s^t \langle Fu_1-Fu_2, 
   u_1-u_2\rangle d\tau~~~~~~~~~~~~~~~~~~~~~~~~~~~ }\\[0.3cm]
 \displaystyle{~~~~~\le \frac 12|u_1(s)-u_2(s)|^2_H+ \int_s^t \langle f_1-f_2,
   u_1-u_2 \rangle d\tau +C^*_0(M)C^*_1(\varepsilon),}
     \end{array} \label{(3.7)}
\end{equation} 
where $C^*_0(M)$ and $C^*_1(\varepsilon)$ are the same ones as in Proposition \ref{prop:2.1}.    
\end{cor}

\noindent
\hspace*{-0.05cm}{\bf Proof.} Just as in the proof of Lemma \ref{lem:3.1}, choose 
approximate sequences 
$\{\theta_{i,n}\} \subset \Theta_S,~\{f_{i,n}\} \subset L^2(0,T;H)$ and 
$\{u_{0i,n}\} \subset H$ for $i=1,2$, such that $u_{0i,n} \in
  \overline{K(\theta_{i,n};0)}$,
 $$ \theta_{i,n} \to \theta_i~{\rm in~}\Theta,~ f_{i,n} \to f_i~{\rm in~}
  L^{p'}(0,T;V^*),~u_{0i,n} \to u_{0i}~{\rm in~} H,~i=1,2,~({\rm as~}n \to \infty).$$
Denoting by $u_{i,n}$ the strong solution of
  $$ u'_{i,n}+\partial \psi^t(\theta_{i,n}; u_{i,n}) \ni f_i,~~u_{i,n}(0)=u_{0i,n}.
  $$
Then, for any large $n$ and small $\varepsilon>0$, 
it follows from \eqref{(3.4)} in the proof of Lemma \ref{lem:3.1}
$$ \begin{array}{l}
   \displaystyle{~~~\frac 12 |u_{1,n}(t)-u_{2,n}(t)|^2_H 
     +\int_s^t \langle Fu_{1,n}-Fu_{2,n}, 
     u_{1,n}-u_{2,n}  \rangle d\tau} \\[0.3cm]
   \displaystyle{\le \frac 12 |u_{1,n}(s)-u_{2,n}(s)|^2_H 
   +\int_s^t (f_{1,n}-f_{2,n}, 
   u_{1,n}-u_{2,n} )_H d\tau + C^*_0(M) C^*_1(\varepsilon),}\\[0.3cm]~~
   \displaystyle{\hspace{5cm}  \forall s,~ t \in [0,T],~s\le t,
     ~~\forall n,~m \geq N_\varepsilon,}
    \end{array} $$ 
By Lemma \ref{lem:3.1}, $\{u_{i,n}\},i=1,2,$ converges in 
$C([0,T];H)\cap L^p(0,T;V)$ to 
solutions $u_i$ of $L_{u_{0i}}(\theta_i;u_i)+Fu_i\ni f_i$ as $n \to \infty$. 
Hence, letting $ n \to \infty$ in the above
inequality, we see that \eqref{(3.7)} holds.\hfill $\Box$\vspace{0.5cm}

\noindent
The following corollary is immediately obtained by letting 
$\varepsilon \downarrow 0$ in \eqref{(3.7)} with $f_1=f_2$ and $u_{10}=u_{20}$.

\begin{cor}
\label{cor:3.2}
 {For every $\theta \in \Theta_W$, 
$u_0 \in \overline{K(\theta;0)}$ and $f \in L^{p'}(0,T;V^*)$, the solution
$u$ of $L_{u_0}(\theta;u)+Fu \ni f$ is unique.}\vspace{0.5cm}
\end{cor}
%


\begin{thm}
\label{th:3.1} 
Assume $(\Theta_S)$ and
(A1)-(A3). Then we have:
\begin{description}
\item{(a)} Let $\theta \in \Theta_W$. Then for any 
$u_0 \in \overline{K(\theta;0)}$, 
$L_{u_0}(\theta;\cdot)$ is a maximal monotone operator from 
$D(L_{u_0}(\theta;\cdot)) \subset L^p(0,T;V)$ into 
$L^{p'}(0,T;V^*)$ and $D(L_{u_0}(\theta;\cdot)) \subset \{w \in C([0,T];H)
\cap L^p(0,T;V)~|~w(0)=u_0,~w(t) \in K(\theta;t)~{\it for~a.e.~}t \in (0,T)\}$.
\item{(b)} Let $\theta \in \Theta_W$, $f,~\bar f \in L^{p'}(0,T;V^*)$,
$u_0,~\bar u_0 \in \overline{K(\theta;0)}$ and $f \in L_{u_0}(\theta; u),
~\bar f \in L_{\bar u_0}(\theta;\bar u)$. Then, for any $s,~t \in [0,T],
~s\le t$,
\begin{equation}
\frac 12|u(t)-\bar u(t)|^2_H \le \frac 12|u(s)-\bar u(s)|^2_H
          +\int_s^t\langle f(\tau)-\bar f(\tau), u(\tau)-\bar u(\tau) 
    \rangle d\tau. \label{(3.8)}
\end{equation} 
\item{(c)} Let $u_0 \in \overline{K(\theta;0)}$ and $f \in L_{u_0}(\theta;u)$.
Then, for any $s, ~t \in [0,T],~s \le t$, 
\begin{equation}
\int_s^t \langle \eta'-f, u-\eta \rangle d\tau
     +\frac 12|u(t)-\eta(t)|^2_H \le \frac 12|u(s)-\eta(s)|^2_H,
   ~~\forall \eta \in {\cal K}_0(\theta). \label{(3.9)}
\end{equation} 
\end{description} 
\end{thm}
\noindent
{\bf Proof.} First we prove (a) and (b). Let $\theta \in \Theta_W$, $u_0 \in
\overline {K(\theta;0)}$ and $f,~\bar f \in L^{p'}(0,T;V^*)$. Assume that
$f \in L_{u_0}(\theta; u)$ and $\bar f \in L_{u_0}(\theta; \bar u)$. Then,
since $f+Fu \in L_{u_0}(\theta; u)+Fu$ and $\bar f+F\bar u 
\in L_{u_0}(\theta;\bar u)+F\bar u$, it follows from \eqref{(3.7)} by letting 
$\varepsilon \to 0$ that  
 \begin{eqnarray*}
  && \frac 12|u(t)-\bar u(t)|^2_H 
    +\int_s^t\langle Fu(\tau)- F\bar u(\tau), u(\tau)-\bar u(\tau) 
    \rangle d\tau \\
  &\le& \frac 12|u(s)-\bar u(s)|^2_H
      +\int_s^t\langle f+Fu(\tau)-\bar f-F\bar u(\tau), u(\tau)-\bar u(\tau) 
    \rangle d\tau, 
  \end{eqnarray*}
which is just \eqref{(3.8)}. By \eqref{(3.8)} with $s=0$,
  $$ \frac 12|u(t)-\bar u(t)|^2_H \le \int_0^t \langle f- \bar f, u - \bar u
       \rangle d\tau,~~\forall t \in [0,T]. $$
This implies that $L_{u_0}(\theta;\cdot)$ is strictly monotone from 
$L^p(0,T;V)$ into $L^{p'}(0,T;V^*)$. Moreover, by Lemma \ref{lem:3.1}, the range of
$L_{u_0}(\theta;\cdot)+F$ is the whole of $L^{p'}(0,T;V^*)$, so that 
$L_{u_0}(\theta;\cdot)$ is maximal monotone from $L^p(0,T;V)$ into 
$L^{p'}(0,T;V^*)$ and $D(L_{u_0}(\theta;\cdot)) \subset \{w\in L^p(0,T;V)\cap
 C([0,T];H)~|~
w(0)= u_0,~w(t) \in K(\theta;t)~{\rm for~a.e.~}t 
\in (0,T)\}$. Thus (a) and (b) are obtained.

Next we show (c). Assume that $f \in L_{u_0}(\theta;u)$. Putting $g:=f+Fu$,
we observe that $g \in L_{u_0}(\theta;u)+Fu$. By \eqref{(3.6)} in the proof
of Lemma \ref{lem:3.1}, we have
  $$ \int_s^t \langle \eta' -g+Fu, u- \eta\rangle d\tau 
    +\frac 12|u(t)-\eta(t)|^2_H
  \le \frac 12 |u(s)-\eta(s)|^2_H, ~~\forall \eta \in {\cal K}_0(\theta).$$
Since $g-Fu=f$, we obtain \eqref{(3.9)}.  \hfill $\Box$ \vspace{0.5cm}

Another important property of $L_{u_0}(\theta;\cdot)$ is stated in the following
theorem. 

\begin{thm}\label{th:3.2}
Assume $(\Theta_S)$ and (A1)-(A3). Let $\{\theta_n\}$ be a sequence in $\Theta_W$ such 
that $\theta_n \to \theta$ in $\Theta$ and let $\{u_{0n}\}$ be a sequence in 
$H$ and $u_0 \in \overline {K(\theta:0)}$ such that
$u_{0n} \in \overline {K(\theta_n:0)}$ for all $n$ and $u_{0n} \to u_0$ in 
$H$ (as $n \to \infty$). Then $\{L_{u_{0n}}(\theta_n;\cdot)\}$ converges to 
$L_{u_0}(\theta;\cdot)$
in the graph sense; namely, if $g \in L_{u_0}(\theta; u)$, then
there are sequences $\{g_n\}$ and $\{u_n\}$ such that 
\begin{equation}
g_n \in L_{u_{0n}}(\theta_n: u_n),~\forall n,~~g_n \to g ~{\it in~}
   L^{p'}(0,T;V^*),~~ u_n \to u~{\it in~} L^p(0,T;V). \label{(3.10)}
\end{equation}
\end{thm}

\noindent
{\bf Proof.} Assume that $g \in L_{u_0}(\theta;u)$. We see that 
$\tilde g:=g+Fu \in L_{u_0}(\theta;u)+Fu$. Take a sequence $\{\tilde g_n\}$
in $L^2(0,T;H)$ so that $\tilde g_n \to g+Fu$ in $L^{p'}(0,T;V^*)$, and 
consider
the sequence $\{u_n\}$ of strong solutions to 
  $$ u'_n(t)+ \partial \psi^t(\theta_n;u_n(t)) \ni \tilde g_n(t)~~{\rm in~}H,
     ~{\rm a.e.~}t \in (0,T),~~u_n(0)=u_{0n},$$
or equivalently
  $$\tilde g_n \in L_{u_{0n}}(\theta_n;u_n) +Fu_n~~{\rm in~}L^{p'}(0,T;V^*).$$
By virtue of Lemma \ref{lem:3.1} and its proof, we observe that $u_n$ converges to the
solution $u$ of $L_{u_0}(\theta; u)+Fu \ni \tilde g$ in $L^{p'}(0,T;V^*)$
in the sense that 
$u_n \to u$ in $C([0,T];H)\cap L^p(0,T;V)$.
Since $\tilde g_n -Fu_n \to \tilde g -Fu=g$ in $L^{p'}(0,T;V^*)$ and the 
sequence $\{u_n, g_n\}$
with $g_n:=\tilde g_n -Fu_n$ satisfies \eqref{(3.10)}. \hfill $\Box$

\section{Parabolic variational inclusions}
\label{sec:pvi}

We begin with the precise assumptions on (multivalued) semimonotone operator 
$A=A(v; w)$. Let ${\cal V}_A$ be a closed convex subset of $L^p(0,T;V)$
such that
\begin{equation}
{\cal K}(\Theta_W):=\bigcup_{\theta \in \Theta_W} {\cal K}(\theta)
    \subset {\cal V}_A.\label{(4.1)}
\end{equation}

Let $A:=A(v;u)$ be a mapping from ${\cal V}_A\times L^p(0,T;V)$ into 
$L^{p'}(0,T;V^*)$,
satisfying 
\begin{description}
\item{(B0)} if $v_1,~v_2 \in {\cal V}_A$ and $w_1,~w_2 \in L^p(0,T;V)$ 
such that $v_1=v_2$ in $V$ and $w_1=w_2$ in $V$ a.e. on $(0,t),~t \in [0,T]$, 
then $[A(v_1; w_1)](\tau)=[A(v_2;w_2)](\tau)$ in $V^*$
for a.e. $\tau \in (0,t)$.

\item{(B1)} (Boundedness) There are positive constants $a_1,a_2$ such that
  $$\sup_{\alpha^*\in A(v,w)}|\alpha^*|^{p-1}_{L^{p'}(0,T;V^*)} \le a_1|w|^{p-1}_{L^p(0,T;V)}+a_2,~~\forall
     v \in {\cal V}_A,~\forall w\in L^p(0,T;V).$$ 
\item{(B2)} (Coerciveness) There are positive constants $a_3,~a_4$ such that
 $$ \int_0^T\langle \alpha^*, w\rangle dt \geq a_3|w|^p_{L^p(0,T;V)}-a_4,~~
     \forall v \in {\cal V}_A,~\forall \alpha^* \in A(v,w).$$
\item{(B3)} (Semimonotonicity) For each $v \in {\cal V}_A$,
$w \to A(v; w)$ is a (multivalued) maximal monotone mapping from 
$D(A(v;\cdot))=L^p(0,T;V)$ into $L^{p'}(0,T;V^*)$. Moreover, for any sequence
$\{v_n\} \subset {\cal V}_A$ with $v_n \to v$ in $L^p(0,T;H)$ and weakly in
$L^p(0,T;V)$ (as $n \to \infty)$, the maximal monotone mapping $A(v_n;\cdot)$ 
converges to
$A(v;\cdot)$ in the graph sense, namely for any
$w \in L^p(0,T;V)$ and any $\alpha^* \in A(v;w)$ there exist sequences
$\{w_n\} \subset L^p(0,T;V)$ and $\{\alpha^*_n\}$ with 
$\alpha^*_n \in A(v_n;w_n)$ such that
   $$ w_n \to w~~{\rm in~}L^p(0,T;V),~~\alpha^*_n \to \alpha^*~~{\rm in~}
      L^{p'}(0,T;V^*). $$
\end{description}

For simplicity, we use the following notation:
  $$ \langle\langle g, w \rangle\rangle =\int_0^T \langle g(t), w(t)\rangle
      dt,~~\forall w \in L^p(0,T;V),~~\forall g \in L^{p'}(0,T;V^*);$$
namely $\langle\langle\cdot,\cdot \rangle\rangle$ means the duality between
$L^{p'}(0,T;V^*)$ and $L^p(0,T:V)$.

For each $\theta \in \Theta_W,~v \in {\cal V}_A,~f \in L^{p'}(0,T;V^*)$ 
and $u_0 \in \overline{K(\theta;0)}$ 
we consider a nonlinear inclusion of the form:
\begin{equation} L_{u_0}(\theta;u) + A(v;u) \ni f~~{\rm in~}L^{p'}(0,T;V^*), \label{(4.2)}
\end{equation}
more precisely, there are $u \in D(L_{u_0}(\theta;\cdot))$, 
$\ell^* \in L_{u_0}(\theta;u)$ and $\alpha^* \in A(v;u)$ such that
   $$ \ell^*(t)+\alpha^*(t)=f(t)~~{\rm in~}V^*,~{\rm a.e.~}t \in (0,T),$$
which is written as
  $$ \begin{array}{l}
  \displaystyle{u\in {\cal K}(\theta),
    ~~\alpha^*\in A(v;u),}\\[0.3cm]
    \displaystyle{ \langle\langle \eta'-f+\alpha^*, u-\eta \rangle\rangle 
    \le \frac 12 |u_0-\eta(0)|^2_H,~~ \forall \eta \in {\cal K}_0(\theta).}
     \end{array}  $$

We now prove:

\begin{thm}
\label{th:4.1}
Assume $(\Theta_S)$, 
(A1)-(A3) and (B0)-(B3).
Then, for each $v \in {\cal V}_A,~\theta \in \Theta_W,~
 f \in L^{p'}(0,T;V^*)$ and 
$u_0 \in \overline{K(\theta;0)}$, there is a unique solution
$u$ of \eqref{(4.2)}. 
\end{thm}

\noindent
{\bf Proof.} By (a) of Theorem \ref{th:3.1}, $L_{u_0}(\theta;\cdot)$ is a
maximal monotone operator from $L^p(0,T;V)$ into $L^{p'}(0,T;V^*)$. Also,
$A(v;\cdot)$ is everywhere defined on $L^p(0,T;V)$, coercive and maximal
monotone from $L^p(0,T;V)$ into $L^{p'}(0,T;V^*)$ by condition (B0)-(B3).
Therefore, it follows from the general theory on monotone operators 
(cf. [6]) that
the range of the sum $L_{u_0}(\theta;\cdot)+A(v;\cdot)$ is the whole 
of $L^{p'}(0,T;V^*)$; in other words, for any $f \in L^{p'}(0,T;V^*)$ 
problem \eqref{(4.2)}
has a solution $u$, which is unique by the strict monotonicity of 
$L_{u_0}(\theta;\cdot)$. \hfill $\Box$ \vspace{0.5cm}

As to the continuous dependence of the solution $u$ of \eqref{(4.2)} upon the 
parameters $\theta$ and $v$, we have:

\begin{thm}\label{th:4.2}
Suppose that $(\Theta_S)$,
 (A1)-(A3) and (B0)-(B3) are fulfilled. Let $f \in L^{p'}(0,T;V^*)$, 
$\theta \in \Theta_W$, $v \in {\cal V}_A$ and 
$u_0 \in \overline{K(\theta;0)}$.
Assume that $\{\theta_n\}\subset \Theta_W$, $\{v_n\} \subset {\cal V}_A$ and 
$\{u_{0n}\}\subset H$ such that $u_{0n} \in \overline{K(\theta_n;0)}$ for all 
$n$, $\{v_n\}$ is bounded in $L^p(0,T;V)$ and
   $$ u_{0n}\to u_0~{\it in~}H,~~
     \theta_n \to \theta~{\it in~}\Theta,~ v_n \to v~{\it in~}L^p(0,T:H) 
    ~({\it as~}n \to \infty). $$
Then,
the sequence $\{u_n\}$ of solutions of \eqref{(4.2)} with $\theta=\theta_n,~v=v_n$ 
and $u_0=u_{0n}$ converges
to the solution $u$ of \eqref{(4.2)} in $C([0,T];H)$ and weakly in $L^p(0,T;V)$.
\end{thm}

By virtue of Theorem \ref{th:4.1}, for each $n$ there is a
unique solution $u_n $ of 
 \begin{equation}  L_{u_{0n}}(\theta_n;u_n)+A(v_n;u_n) \ni f ~~{\rm in~}L^{p'}(0,T;V^*),
    \label{(4.3)}
    \end{equation}
namely
  $$ \ell^*_n \in L_{u_0}(\theta_n;u_n),~\alpha^*_n \in A(v_n;u_n),~~
    \ell^*_n(t)+\alpha^*_n(t) =f(t)~~{\rm in~}V^*,~{\rm a.e.~}t \in(0,T).
   $$

\begin{lem}
\label{lem:4.1} 
The sequence $\{u_n\}$ of solutions to \eqref{(4.3)} is bounded
in $L^p(0,T;V)$ and $C([0,T]; H)$; in fact we have:
\begin{equation} |u_n|^p_{L^p(0,T;V)}+ |u_n|^2_{C([0,T];H)}
   \le N_0 \left\{|f|^{p'}_{L^{p'}(0,T;V^*)}+|u_0|^2_H + 1\right\},~~
   \forall n,
     \label{(4.4)}
    \end{equation}
where $N_0$ is a positive constant independent of $n$. 
\end{lem}

\noindent
{\bf Proof.} By Lemma \ref{lem:2.1} there is a sequence $\{\eta_n\}$ with $\eta_n \in
{\cal K}_0(\theta_n)$ with a positive constant $N'_0$ such that
  $$ |\eta_n|_{L^p(0,T;V)} + |\eta'_n|_{L^{p'}([0,T];V^*)} +|\eta_n|_{C(0,T;H)}
    \le N'_0.
  $$
For each $n$ we have by (c) of Theorem \ref{th:3.1}
  $$ \int_0^t \langle \eta'_n, u_n-\eta_n \rangle d\tau
  +\int_0^t \langle \alpha^*_n, u_n-\eta_n\rangle 
     +\frac 12|u_n(t)-\eta_n(t)|^2_H$$
 $$\le \int_0^t \langle f, u_n-\eta_n \rangle d\tau 
    +\frac 12|u_{0n}-\eta_n(0)|^2_H,~~\forall t \in [0,T].$$
From the above inequality with (B1) and (B2) we obtain the following estimate:
  $$ \frac {a_3}2 \int_0^t |u_n|^p_Vd\tau + \frac 14|u_n(t)|^2_H 
   \le N''_0 \left\{\int_0^t |f|^{p'}_{V^*}d\tau+|u_0|^2_H +1\right\},
    ~~\forall t \in [0,T],$$
where $a_3$ is the same constant as in (B2), $N''_0$ is a positive constant 
independent of
$t \in [0,T],~\eta_n$ and $f$; actually it depends only on $N'_0$. Hence we 
have \eqref{(4.4)}.
\hfill $\Box$ \vspace{0.5cm}

On account of Lemma \ref{lem:4.1}, we can find a subsequence $\{n_k\}$
of $\{n\}$ such that
  $$ u_{n_k} \to u ~{\rm weakly~in~}L^p(0,T;V),~~\alpha^*_{n_k} \to 
   \alpha^*~{\rm weakly~in~}L^{p'}(0,T;V^*),$$
 $$\ell^*_{n_k}=f-\alpha^*_{n_k} \to f-\alpha^*=:\ell^*~{\rm weakly~in~}
  L^{p'}(0,T;V^*), $$
as $k\to \infty$.

\begin{lem}
\label{lem:4.2}
$\displaystyle{~\liminf_{k\to \infty}
    \langle\langle \alpha^*_{n_k}, u_{n_k}-u \rangle \rangle \geq 0~~
  and~~\limsup_{k\to \infty} \langle\langle \ell^*_{n_k}, u_{n_k}-u
    \rangle\rangle \le 0 .}$  
\end{lem}

\noindent
{\bf Proof.} We first observe from \eqref{(4.3)} that
\begin{equation}  \langle\langle f, u_{n_k}-w \rangle\rangle 
   =\langle\langle \ell^*_{n_k}, u_{n_k}-w \rangle\rangle
    + \langle\langle \alpha^*_{n_k}, u_{n_k}-w \rangle\rangle,~~\forall w
    \in L^p(0,T;V).\label{(4.5)}
    \end{equation}
By assumption (B3), for any $\tilde \alpha^* \in A(v;u)$ there are sequences 
$\{w_k\}$ in $L^p(0,T;V)$ and
$\{\tilde\alpha^*_k\}$ with $\tilde\alpha^*_k \in A(v_{n_k};w_k)$ such that
  $$ w_k \to u~~{\rm in~}L^p(0,T;V),~~\tilde\alpha^*_k \to \tilde\alpha^*~{\rm in~}
     L^{p'}(0,T;V^*).$$
Now, taking $w_k$ as $w$ in \eqref{(4.5)} and passing to the limit as $k\to \infty$,
we have
\begin{equation} 0 = \lim_{k\to \infty}\langle\langle f, u_{n_k}-w_k \rangle\rangle
      \geq \limsup_{k \to \infty}\langle\langle \ell^*_{n_k}, u_{n_k}-w_k
           \rangle\rangle +\liminf_{k \to \infty}\langle\langle \alpha^*_{n_k},              u_{n_k}-w_k \rangle\rangle. \label{(4.6)}
    \end{equation}
Here, since $\alpha^*_{n_k} \in A(v_{n_k};u_{n_k})$ and $\tilde\alpha^*_k \in 
A(v_{n_k}; w_k)$, we note from the monotonicity of $A(v_{n_k};\cdot)$ that
 \begin{eqnarray*}
 \langle\langle \alpha^*_{n_k}, u_{n_k}-w_k \rangle\rangle
 &=&\langle\langle \alpha^*_{n_k}-\tilde\alpha^*_k, u_{n_k}-w_k \rangle\rangle
   +\langle\langle \tilde \alpha^*_k, u_{n_k}-w_k \rangle\rangle \\
 &\geq& \langle\langle \tilde \alpha^*_k, u_{n_k}-w_k \rangle\rangle. 
 \end{eqnarray*}
 Therefore,
  $$ \liminf_{k\to \infty}\langle\langle \alpha^*_{n_k}, u_{n_k}-w_k 
     \rangle\rangle \geq \lim_{k\to \infty}\langle\langle \tilde\alpha^*_k,
     u_{n_k}-w_k \rangle\rangle =\langle\langle \tilde \alpha^*, u-u
      \rangle \rangle =0,$$
so that
  \begin{eqnarray*}
   \liminf_{k\to \infty}\langle\langle \alpha^*_{n_k}, u_{n_k}-u 
     \rangle\rangle &=&\liminf_{k\to \infty}\langle\langle \alpha^*_{n_k}, 
    u_{n_k}-w_k +w_k-u \rangle\rangle \\
    &=&\liminf_{k\to \infty}\langle\langle \alpha^*_{n_k}, u_{n_k}-w_k
     \rangle\rangle 
   + \lim_{k\to \infty}\langle\langle \alpha^*_{n_k}, w_k-u\rangle\rangle\\
   &=&\liminf_{k\to \infty}\langle\langle \alpha^*_{n_k}, u_{n_k}-w_k
     \rangle\rangle \geq 0
  \end{eqnarray*}
and another inequality of the lemma follows similarly from \eqref{(4.6)}. \hfill 
$\Box$ \vspace{0.3cm}

\noindent
{\bf Proof of Theorem \ref{th:4.2}:} First we shall show that 
\begin{equation} \lim_{k\to \infty}\langle\langle \alpha^*_{n_k},
      u_{n_k}\rangle\rangle =\langle\langle \alpha^*, u\rangle\rangle~~
   {\rm and~~} \alpha^* \in A(v;u). \label{(4.7)} 
    \end{equation}
Let $\eta $ be any function in ${\cal K}_0(\theta)$ and put
$\eta_{n,\varepsilon}(t):={\cal F}_{\theta\theta_n, \varepsilon}(t)
\eta(t)$; note from Lemma \ref{lem:2.1} that $\eta_{n,\varepsilon} \to \eta$ 
in $L^p(0,T;V),$
$\eta'_{n,\varepsilon} \to \eta'$ in $L^{p'}(0,T;V^*)$ and 
$\eta_{n,\varepsilon} \to \eta$ in $C([0,T];H)$ as well,
when $n \to \infty$ and $\varepsilon \downarrow 0$.
Since $f-\alpha^*_n \in L_{u_{0n}}(\theta_n; u_n)$, it follows that
\begin{equation}  \langle\langle \eta'_{n,\varepsilon} -f+\alpha^*_n, u_n
    -\eta_{n,\varepsilon} \rangle\rangle
  \le \frac 12|u_{0n}-\eta_{n,\varepsilon}(0)|^2_H,
    \label{(4.8)}
    \end{equation}
Now, let $n=n_k \to \infty$ and $\varepsilon \downarrow 0$ in \eqref{(4.8)}
and note $\liminf_{n\to \infty}\langle\langle \alpha^*_n,u_n\rangle\rangle
\geq \langle\langle \alpha^*, u \rangle\rangle$ by Lemma \ref{lem:4.2}, to see that
  $$ \langle\langle \eta' -f+\alpha^*, u-\eta \rangle\rangle
  \le \frac 12|u_0-\eta(0)|^2_H,~~ \forall \eta \in {\cal K}_0(\theta),
   $$
which implies by definition that $f-\alpha^* \in L_{u_0}(\theta;u)$.

Next we show that $\alpha^* \in A(v;u)$. In fact, since $\ell^*:=f-\alpha^*
\in L_{u_0}(\theta;u)$, it follows from Theorem \ref{th:3.2} that there exists
a sequence $\{\tilde u_n, \tilde \ell^*_n\}$ such that 
 \begin{equation} \tilde\ell^*_n \in 
    L_{u_{0n}}(\theta_n;\tilde u_n),~\tilde u_n \to u~{\rm in~}L^p(0,T;V),~
    \tilde\ell^*_n \to \ell^*=f-\alpha^* ~{\rm in~}L^{p'}(0,T;V^*).
   \label{(4.9)}
    \end{equation}
Using this sequence, we see that
  \begin{eqnarray*}  
  \liminf_{k \to \infty}\langle\langle \ell^*_{n_k}, u_{n_k}-u\rangle\rangle
     &=&\liminf_{k \to \infty}\langle\langle \ell^*_{n_k}, u_{n_k}
         -\tilde u_{n_k}+\tilde u_{n_k} -u \rangle\rangle\\
     &\geq& \liminf_{k \to \infty}\langle\langle \tilde \ell^*_{n_k}, 
          u_{n_k}-\tilde u_{n_k}\rangle\rangle
    +\liminf_{k \to \infty}\langle\langle \ell^*_{n_k}, 
     \tilde u_{n_k} -u \rangle\rangle \\
    &\geq& \lim_{k \to \infty}\langle\langle \tilde \ell^*_{n_k}, 
          u_{n_k}-\tilde u_{n_k}\rangle\rangle
    +\lim_{k \to \infty}\langle\langle \ell^*_{n_k}, 
     \tilde u_{n_k} -u \rangle\rangle \\ 
    &=& 0 
   \end{eqnarray*}
Therefore, together with 
$\limsup_{k\to \infty}\langle\langle \ell^*_{n_k}, u_{n_k}
-u\rangle\rangle\le 0$ in Lemma \ref{lem:4.2},
we have that 
  $$\lim_{k \to \infty}\langle\langle \ell^*_{n_k}, u_{n_k}-u\rangle\rangle 
    =0, ~~{\rm i.e.~}\lim_{k \to \infty}\langle\langle \ell^*_{n_k}, u_{n_k}
    \rangle\rangle =\langle\langle \ell^*,u \rangle\rangle. $$
Since $\ell^*_{n_k}=f-\alpha^*_{n_k}$ and $\alpha^*_{n_k} \to \alpha^*$
weakly in $L^{p'}(0,T;V^*)$, it results from the above equality that
\begin{equation}  \lim_{k \to \infty}\langle\langle \alpha^*_{n_k}, u_{n_k}\rangle\rangle
    = \langle\langle  \alpha^*, u \rangle\rangle, \label{(4.10)}
    \end{equation}
namely \eqref{(4.7)} holds.

We are now in a position to show $\alpha^* \in A(v;u)$. 
Let $w$ be any 
element in $V$ and $\alpha^*_w$ be any element of $A(v;w)$. By (B3), choose
a sequence $\{\alpha^*_{w,k}, w_k\}$ so that $\alpha^*_{w,k} 
\in A(v_{n_k};w_k)$, $w_k \to w$ in $L^p(0,T;V)$ and $\alpha^*_{w,k} \to
\alpha^*_w$ in $L^{p'}(0,T;V^*)$ as $k\to \infty$.
 Then, by the
monotonicity of $A(v_{n_k};\cdot)$, we have that
  $$ \langle\langle \alpha^*_{n_k}-\alpha^*_{w,k}, u_{n_k}-w_k \rangle \rangle
     \geq 0.$$
Passing to the limit in $k\to \infty$, we obtain by \eqref{(4.10)} 
   $$ \langle\langle \alpha^*-\alpha^*_w, u-w \rangle \rangle \geq 0,
   ~~\forall \alpha^*_w \in A(v;w),~\forall w \in V,$$
which implies that $\alpha^* \in A(v;u)$ by the maximal monotonicity of 
$A(v;\cdot): V\to V^*$.

Finally we show that $u_{n_k} \to u$ in $C([0,T];H)$ as $k \to \infty$. From
(b) of Theorem \ref{th:3.1} it follows that for all $t \in [0,T]$
\begin{equation} \frac 12 |u_{n_k}(t)-u(t)|^2_H 
        +\int_0^t \langle \alpha^*_{n_k}-\alpha^*, u_{n_k}-u\rangle d\tau
    \le \frac 12 |u_{0n_k}-u_0|^2_H
    +C^*_0(M)C^*_1(\varepsilon), \label{(4.11)}
    \end{equation}
for all large $k$ and any small $\varepsilon >0$. 

Here we show by the same 
idea as getting \eqref{(4.10)} under \eqref{(4.9)} that
\begin{equation}  \lim_{k\to \infty}\int_0^t \langle \alpha^*_{n_k}-\alpha^*, 
    u_{n_k}-u\rangle d\tau = 0.
   \label{(4.12)}
    \end{equation}
Indeed, take sequences $\{\tilde u_n\}$ and $\{\tilde \alpha^*_n\}$
so that $\tilde u_n \to u$ in $L^p(0,T;V)$, $\tilde\alpha^*_n \to \alpha^*$ 
in $L^{p'}(0,T;V^*)$ and $\tilde \alpha^*_n \in A(v_n; \tilde u_n)$ for all
$n$. Then, by putting 
  $$\bar u_{n_k}:=\left\{
     \begin{array}{ll}
      \tilde u_{n_k}&~~~{\rm on~}[0,t),\\
      u_{n_k}&~~~{\rm on~}[t,T], 
     \end{array} \right.\qquad
    \bar \alpha^*_{n_k}(t):=\left\{
     \begin{array}{ll}
       \tilde\alpha^*_{n_k}&~~~{\rm on~}[0,t),\\
       \alpha^*_{n_k}&~~~{\rm on~} [t,T],
     \end{array} \right. $$
we see from condition (B0) that 
$\bar\alpha^*_{n_k}\in A(v_{n_k};\bar u_{n_k})$,
and clearly $\bar u_{n_k} \to u$ in $L^p(0,T;V)$ and
$\bar \alpha^*_{n_k} \to \alpha^*$ in $L^{p'}(0,T;V^*)$ 
 as $k \to \infty$. Hence we can obtain 
 \begin{eqnarray*} 
 &&\liminf_{k\to \infty}\int_0^t \langle \alpha^*_{n_k}-\alpha^*, 
   u_{n_k}-u\rangle d\tau \\
 &=&\liminf_{k\to \infty}\int_0^t \langle \alpha^*_{n_k}-\tilde \alpha^*_{n_k}
   +\tilde \alpha^*_{n_k}-\alpha^*,u_{n_k}-\tilde u_{n_k}+\tilde u_{n_k}
    -u\rangle d\tau \\
 &=&\liminf_{k\to \infty}\int_0^t \langle 
   \alpha^*_{n_k}-\tilde \alpha^*_{n_k}, u_{n_k}-\tilde u_{n_k} \rangle d\tau\\
& =&\liminf_{k\to \infty}\int_0^T \langle \alpha^*_{n_k}-\bar \alpha^*_{n_k}, 
   u_{n_k}- \bar u_{n_k}\rangle d\tau \geq 0;
 \end{eqnarray*}
the last equality follows from (B3). As $k \to \infty$ and 
$\varepsilon \downarrow 0$ in \eqref{(4.11)}, \eqref{(4.12)} is obtained and $u_{n_k} \to u$ 
in $C([0,T]; H)$.
\hfill $\Box$ 

\section{Parabolic quasi-variational inequalities}
\label{sec:pqvi}

In this section we give the formulation of a class of parabolic
quasi-variational inequalities and an abstract existence result.

Given an initial datum $u_0 \in H$, we introduce a feedback system 
$\Lambda_{u_0}$, which is an operator from
${\cal V}_A$ (cf. \eqref{(4.1)}) into $\Theta_W$ satisfying the following conditions:
\begin{description}
\item{($\Lambda1$)} $\Lambda_{u_0}$ maps ${\cal V}_A$ into 
$\Theta_W(u_0):=\{\theta \in \Theta_W~|~u_0 \in \overline{K(\theta;0)}\}$.
\item{($\Lambda2$)} If $\{w_n\} \subset {\cal V}_A$ and is bounded in 
$L^p(0,T;V)$ 
and $w_n \to w$ in $L^p(0,T; H)$, then $\Lambda_{u_0} w_n \to \Lambda_{u_0}w$
in $\Theta$ as $n \to \infty$. 
\end{description}

\begin{defn}
\label{def:5.1}
 Given $u_0 \in H$ and $f \in L^{p'}(0,T;V^*)$, we denote
by $QVI(\Lambda_{u_0}; f, u_0)$ the problem to find a pair $\{\theta,u\} 
\in \Theta_W(u_0) \times L^p(0,T;V)$ such that
\begin{equation}  \left\{ \begin{array}{l}
  \displaystyle{ L_{u_0}(\theta;u)+A(u;u) \ni f~~{\rm in~}L^{p'}(0,T;V^*),}
     \\[0.3cm]
  \displaystyle{ \theta=\Lambda_{u_0} u~{\rm in~}\Theta_W.}
 \end{array} \right. \label{(5.1)}
    \end{equation}
\end{defn}

We need an additional technical set-up in order to establish the solvability of
$QVI(\Lambda_{u_0};$
$f,u_0)$. Let $W$ be a reflexive and {\it separable} Banach 
space 
which is densely and continuously embedded in $V$; in this case, since 
$V^*\subset W^*$,
 \begin{equation}  V \subset H \subset W^*~{\rm with ~compact~ embeddings}. \label{(5.2)}
    \end{equation}
Moreover, we suppose that there is a positive number $\delta_0$ such that
\begin{equation}   \delta_0 B_W(0) \subset K(\theta;t),~~\forall \theta \in \Theta_W(u_0),
   ~\forall  t \in [0,T], \label{(5.3)} 
    \end{equation}
where $B_W(0)$ is the unit closed ball around the origin in $W$.

Under the additional conditions \eqref{(5.2)}, \eqref{(5.3)} and a given constant $M^*>0$, 
we consider the set 
$Z(\delta_0,M^*, u_0)$ in $L^p(0,T;V)\cap L^\infty(0,T;H)$ given as:
\begin{equation}  \left\{u \left|~
   \begin{array}{l}
  \displaystyle{|u|_{L^p(0,T;V)} \le M^*,~|u|_{L^\infty(0,T;H)} \le M^*,}
       \\[0.2cm]
  \displaystyle{ \exists g \in L^{p'}(0,T;V^*)~{\rm such~that~}}\\[0.2cm]
  \displaystyle{ ~~~\langle\langle g,u \rangle\rangle \le M^*,
  ~~|g|_{L^1(0,T;W^*)} \le M^*,} \\[0.2cm]
   \displaystyle{~~~\langle\langle \eta'-g, u-\eta \rangle\rangle 
    \le \frac 12|u_0-\eta(0)|^2_H,~~\forall \eta \in L^p(0,T;V)}\\[0.2cm]
  \displaystyle{~~~~~~{\rm with~}\eta' \in L^{p'}(0,T;V^*)~{\rm and~}\eta(t) 
   \in \delta_0 B_W(0),~   \forall t \in [0,T]} 
   \end{array} \right. \right\} \label{(5.4)} 
    \end{equation}
The compactness lemma stated below is one of important mathematical tools for
the solvability of $QVI(\Lambda_{u_0}; f, u_0)$. 

\begin{lem}\label{lem:5.1}
{\bf \rm ([10; Theorem 4.1]).}
 For any $\delta_0>0$, $M^*>0$ and
$u_0 \in H$ the set $Z(\delta_0,M^*, u_0)$
is relatively compact in $L^p(0,T;H)$ and its convex closure $\overline
{\rm conv}(Z(\delta_0,M^*,u_0))$ (in $L^p(0,T;V)$) is compact in 
$L^p(0,T;H)$.
\end{lem}

Next, let $A(v,u)$ be the same semimonotone operator as in the last section.
Then we note that there is a positive constant $N_0$ 
such that
\begin{equation}  |u|^2_{C([0,T];H)} + |u|^p_{L^p(0,T; V)} \le N_0(|f|^{p'}_{L^{p'}(0,T;V^*)}   + |u_0|^2_H) =:N_1, \label{(5.5)} 
    \end{equation}
for all solutions $u$ of $L_{u_0}(\theta;u)+A(v;u) \ni f$ as long as 
$\theta \in \Theta_W(u_0)$ and $v \in {\cal V}_A$. In fact, by (c) of 
Theorem \ref{th:3.1} we have
 $$ \int_0^t \langle \eta' -f+\alpha^*, u-\eta \rangle d\tau 
    +\frac 12|u(t)-\eta(t)|^2_H \le \frac 12 |u_0-\eta(0)|^2_H,
    ~\forall \eta \in {\cal K}_0(\theta),~\forall t \in [0,T],$$
where $\alpha^* \in A(v;u)$.
By \eqref{(5.3)}, we take $0$ as $\eta \in {\cal K}_0(\theta_n)$ to get
 $$ \frac 12|u(t)|^2_H
  + \int_0^t \langle \alpha^*,u \rangle d\tau 
  \le \int_0^t\langle f, u \rangle d\tau+ \frac 12|u_0|^2_H,
  ~~\forall t \in [0,T].  $$
By using conditions (B1) and (B2), it is easy to derive \eqref{(5.5)} for some constant
$N_1>0$ from this inequality.

\begin{cor}\label{cor:5.1}
 Let $\{\theta_n\} \subset \Theta_W(u_0)$ and
$\{v_n\} \subset {\cal V}_A$ be sequences such that $\{v_n\}$ is bounded
in $L^p(0,T;V)$ and 
  $$ \theta_n \to \theta~~{\it in~}\Theta,~v_n \to v~{\it in~}L^p(0,T;H)~~
    ({\it as~}n \to \infty). $$
Let $\{u_n\}$ be the sequence of solutions $u_n$ of $L_{u_0}(\theta_n;u_n)
+A(v_n;u_n) \ni f$.
Then $\{u_n\}$ is bounded in $C([0,T];H)\cap L^p(0,T;V)$ and is relatively 
compact in $L^p(0,T;H)$.
\end{cor}

We see by (c) of Theorem \ref{th:3.1}, together with \eqref{(5.2)}-\eqref{(5.5)}, and Theorem \ref{th:4.2}
that $u_n \in Z(\delta_0,M^*, u_0)$ for a certain constant $M^*>0$. Hence
this corollary is a direct consequence of Lemma \ref{lem:5.1}. 
\medskip

Now we formulate an existence result for $QVI(\Lambda_{u_0};f,u_0)$.

\begin{thm}\label{th:5.1}
Assume that $(\Theta_S)$, (A1)-(A3), 
(B0)-(B3) and \eqref{(5.2)}-\eqref{(5.3)} are fulfilled. Let $u_0\in H$ be a given initial 
datum and
$\Lambda_{u_0}$ be a feedback system from ${\cal V}_A$ into 
$\Theta_W(u_0) \ne \emptyset$,
satisfying $(\Lambda1)$ and $(\Lambda 2)$.
Then, for each $f \in L^{p'}(0,T;V^*)$, problem 
$QVI(\Lambda_{u_0};f,u_0)$ admits at least one solution $\{\theta,u\}$.
\end{thm}

\noindent
{\bf Proof.} By estimate \eqref{(5.5)}, for a large constant $M^*$ we see that
any solution $u$ of $L_{u_0}(\theta;u)+A(v;u) \ni f$ belongs to
$Z(\delta_0,M^*, u_0)$, as long as $\theta\in \Theta_W(u_0)$ and 
$v \in {\cal V}_A$.

We put ${\cal X}(u_0):=\overline{{\rm conv}}(Z(\delta_0;M^*,u_0))$,
which is compact and convex in $L^p(0,T;H)$ by Lemma \ref{lem:5.1}. 
Now, for each $u \in {\cal X}(u_0)\cap {\cal V}_A$, we denote by $\bar u$ 
a unique solution of
   $$ \theta=\Lambda_{u_0}u,~~L_{u_0}(\theta; \bar u)+A(u;\bar u) \ni 
    f~~{\rm in~}L^{p'}(0,T;V^*), $$
and define a mapping $S: {\cal X}(u_0)\cap {\cal V}_A \to {\cal X}(u_0)\cap
{\cal V}_A$ by $\bar u=Su ~(\in
Z(\delta_0,M^*,u_0) \subset {\cal X}(u_0)\cap {\cal V}_A)$.

We are going to prove that $S$ is continuous in ${\cal X}(u_0)\cap {\cal V}_A$
 with respect to the topology of $L^p(0,T;H)$.
Let $\{u_n\}$ be any sequence in ${\cal X}(u_0)\cap {\cal V}_A$ such that 
$u_n \to u$ in
$L^p(0,T;H)$. Then, by $(\Lambda2)$, $\theta_n :=\Lambda_{u_0}u_n \to 
\Lambda_{u_0}u=:\theta$ in $\Theta$. Therefore, by Corollary \ref{cor:5.1}, the solution
$\bar u_n$ of $L_{u_0}(\theta_n;\bar u_n)+A(u_n;\bar u_n) \ni f$ converges
to the solution $\bar u$ of $L_{u_0}(\theta; \bar u)+A(u;\bar u) \ni f$ in
the sense that $\bar u_n \to \bar u$ in $C([0,T];H)$ and weakly in 
$L^p(0,T;V)$. This shows that $\bar u_n=Su_n \to Su=\bar u$ in $L^p(0,T;H)$.

Now we are in a position to apply the Schauder's fixed-point theorem for $S$ in
${\cal X}(u_0)\cap {\cal V}_A$ in order to find a function 
$u \in {\cal X}(u_0)\cap {\cal V}_A$ such that
$u=Su$. In this case, the pair $\{\theta, u\}$ with $\theta=\Lambda_{u_0}u$
satisfies \eqref{(5.1)} and $u$ is a solution of $QVI(\Lambda_{u_0}; f, u_0)$. 
\hfill $\Box$
\vspace{0.5cm}

In general, the quasi-variational inequality \eqref{(5.1)} has multiple solutions
as the following  simple example shows. 

\begin{ex}
(Non-uniqueness)
We consider the case $H=V={\bf R}$. For a fixed positive
constant $c_0$, put
 $$ X_0:=\{z \in W^{1,2}(0,T)~|~ 0\le z' \le c_0~{\rm a.e.~on~}(0,T),~
 z(0)=1\},$$
which is compact and convex in $L^2(0,T)$. Clearly
  \begin{equation}  \{2- e^{-ct}~|~0\le c\le c_0\} \subset X_0. \label{(5.6)}
    \end{equation}
As the space $\Theta$ of parameters we take $X_0$
with metric $d_\Theta(z_1,z_2)=|z_1-z_2|_{C(\overline Q)}$, and define
for every $z \in \Theta$
  $$ K(z;t):=\{r \in {\bf R} ~|~
     r \geq z(t)-2\},~\forall t \in [0,T].$$
It is easy to check the conditions $(\Theta_S)$ and (A1)-(A3).
By the general results in section \ref{sec:L}, the time-derivative $L_1(z; \cdot)$
is defined as well corresponding to the initial value $1$ and the constraint 
set $K(z;t)$.
Next we formulate $\Lambda_1$, with initial value $1$, as a feedback system 
from $L^2(0,T)$ into $\Theta$ 
as follows: for each $v \in L^2(0,T)$, define $z:=\Lambda_1 v \in \Theta$ by
 $$ |z-v|^2_{L^2(0,T)} = \min_{\zeta \in \Theta}|\zeta-v|^2_{L^2(0,T)};$$
we note from \eqref{(5.6)} that
\begin{equation}  \Lambda_1(2-e^{-ct}) =2-e^{-ct}, ~~0\le \forall c \le c_0. \label{(5.7)}
    \end{equation}

Now, consider the quasi-variational inequality \eqref{(5.1)} with $A(v;u)\equiv 0$
and $f=0$:
  $$ L_1(z;u) \ni 0~~{\rm in~}L^2(0,T),~~z=\Lambda_1 u,$$
which is written as
 \begin{equation}  \int_0^T \eta'(u-\eta)dt \le \frac 12|1-\eta(0)|^2,~~
     \forall \eta \in {\cal K}_0(z),~~z=\Lambda_1 u, \label{(5.8)}
    \end{equation}
where ${\cal K}_0(z)$ is the class of smooth test functions, namely
${\cal K}_0(z)=\{\eta \in W^{1,2}(0,T)~|~\eta(t) \in K(z;t),~\forall t 
\in [0,T]\}$.
For instance, pay attention to the function $u(t):=2-e^{-ct}$ for any constant
$c$ with $0\le c \le c_0$. Then we see that $u(0)=1$ and for any 
$\eta \in {\cal K}_0(u)$
 $$ \int_0^T \eta'(u-\eta) dt \le \frac 12|1-\eta(0)|^2. $$
In fact, since $u'(t)=ce^{-ct} \geq 0$ and $u(t)-\eta(t)\le 0$, it follows that
 \begin{eqnarray*}
  \int_0^T \eta'(u-\eta)dt 
  &=&\int_0^T u'(u-\eta)dt -\frac 12 |u(T)-\eta(T)|^2+\frac 12 |1-\eta(0)|^2\\
  &\le& \frac 12 |1-\eta(0)|^2.                            
 \end{eqnarray*}
Besides, we have by \eqref{(5.7)} that $\Lambda_1 u=2-e^{-ct} =u$.
Thus $u:= 2- e^{-ct}$ satisfies \eqref{(5.8)} with $z=u$ for all $c \in [0,c_0]$, 
which shows that problem \eqref{(5.8)} possesses infinite many solutions.
\end{ex}

\begin{req}\label{rem:5.1}
In Theorem \ref{th:5.1} an existence result of $QVI(\Lambda_{u_0};f,
u_0)$ was established, based on a compactness property of
$L_{u_0}(\theta;\cdot)$ (cf. Lemma \ref{lem:5.1}). Of course, there are a variety of 
existence results for $QVI(\Lambda_{u_0};f,u_0)$ without such a compactness 
property, depending on the choice of feedback system $\Lambda_{u_0}$.
For instance, see Application \ref{app1} of the next section.
\end{req}

\section{Applications}
\label{sec:app}

\subsection{Quasi-variational ordinary differential inequality}
\label{app1}

In the first application we are going to treat a sweeping process with 
quasi-variational structure in two dimensional space ${\bf R}^2$.

Let us consider the case of $H=V=W={\bf R}^2$ and 
 $$ \Theta :=\left\{\theta:=[{\bfa}, \gamma, {\boldsymbol\zeta}]~\left|
    \begin{array}{l}
     {\bfa}\in W^{1,p}(0,T;{\bf R}^2),~|{\bfa}(t)| = 1,~\forall
      t \in [0,T],\\
  \gamma \in C_b({\bf R}^2),~\gamma_* \le \gamma \le \gamma^*
     ~{\rm on~}{\bf R},\\
  {\boldsymbol\zeta} \in C([0,T]; {\bf R}^2)
    \end{array} \right. \right \},$$
where $2\le p <\infty$ and $\gamma_*,~\gamma^*$ are positive constants with
$\gamma_*<\gamma^*$; $C_b({\bf R}^2)$ is the space of all functions
$\gamma$ in $C({\bf R}^2)$ such that $\lim_{|{\bfr}|\to \infty}\gamma({\bfr})$ 
exists. Here the space $\Theta$
is a complete metric space with metric
 $$ d_{\Theta}(\theta,\bar\theta):=|{\bfa}-\bar{\bfa}|
  _{W^{1,p}(0,T;{\bf R}^2)}+|\gamma-\bar\gamma|_{C_b({\bf R}^2)}
   +|{\boldsymbol\zeta}-\bar{\boldsymbol\zeta}|_{C([0,T];{\bf R}^2)},$$
for $\theta:=[{\bfa}, \gamma, {\boldsymbol\zeta}],
  ~\bar\theta:=[\bar{\bfa},\bar\gamma, \bar{\boldsymbol\zeta}] \in \Theta$. 
Also, given $\theta:=[{\bfa},\gamma,{\boldsymbol\zeta}]$ and
$\bar\theta:=[\bar{\bfa},\bar\gamma,\bar{\boldsymbol\zeta}]$, we define the 
rotation $R_{\theta\bar\theta}(t)$ by
\begin{equation}  R_{\theta\bar\theta}(t):= \left( \begin{array}{l}
                     \displaystyle{ \cos \alpha(t)~~-\sin \alpha(t)}\\[0.2cm]
                     \displaystyle{\sin \alpha(t)~~~~~ \cos \alpha(t)}
                                     \end{array} \right ), \label{(6.1)}
    \end{equation}
with the angle $\alpha(t)$ between vectors 
${\bfa}(t):=(a^{(1)}(t),a^{(2)}(t))$ and 
$\bar{\bfa}(t):=(\bar a^{(1)}(t),\bar a^{(2)}(t))$.  
It is easy to see from \eqref{(6.1)} that for any vector 
${\bfz}:=(z^{(1)},z^{(2)}) \in {\bf R}^2$
\begin{equation}  R_{\theta\bar\theta}(t){\bfz}
   =(z^{(1)}\cos \alpha(t)-z^{(2)}\sin \alpha(t), z^{(1)}\sin \alpha(t)
    +z^{(2)}\cos \alpha(t)),\label{(6.2)}
    \end{equation}
with
\begin{equation}  \sin \alpha(t)=a^{(1)}(t)\bar a^{(2)}(t)-a^{(2)}(t)\bar a^{(1)}(t),
  ~~\cos \alpha(t)=a^{(1)}(t)\bar a^{(1)}(t)+a^{(2)}(t)\bar a^{(2)}(t).
   \label{(6.3)}
    \end{equation}                
Also, as function $\sigma_{\theta\bar\theta,\varepsilon}(t)$ we take
 \begin{equation}  \sigma_{\theta\bar\theta,\varepsilon}(t):= \varepsilon \bar{\bfa}(t),
     ~~\forall \varepsilon \in (0,1],~\forall t \in [0,T].
    \label{(6.4)}
    \end{equation}

Now, for each $\theta:=[{\bfa}, \gamma, {\boldsymbol\zeta}] \in \Theta$, we put
\begin{equation} K(\theta;t):=\{{\bfz}\in {\bf R}^2~|~{\bfa}(t)\cdot({\bfz}-{\bfa}(t))=0,~
       |{\bfz}-{\bfa}(t)|\le \gamma({\boldsymbol\zeta}(t))\},
  ~~\forall t \in [0,T].
    \label{(6.5)}
    \end{equation}

\noindent
({\it Verification of $(\Theta_S)$}) Denote by $\Theta_1$ 
the set of all 
parameters $\theta:=[{\bfa},\gamma,{\boldsymbol\zeta}]\in \Theta$ of 
$C^2$-class.
Let $\theta:=[{\bfa},\gamma,{\boldsymbol\zeta}] \in \Theta_1$, and 
$0=T_0<T_1<T_2<\cdots<
T_N:=T$ be a partition of $[0,T]$ such that 
  $$ |\gamma({\boldsymbol\zeta}(s))-\gamma({\boldsymbol\zeta}(t))|
   < \gamma_*,~~\forall s,~t \in
     [T_{k-1},T_k],~k=1,2,\cdot, N. $$
Given ${\bfz} \in K(\theta;s),~s,~t \in [T_{k-1},T_k],~s\le t$, we put
 $$\tilde {\bfz}:=\left(1-\frac 1{\gamma_*}|\gamma({\boldsymbol\zeta}(s))
  -\gamma({\boldsymbol\zeta}(t))|\right)R(s,t)({\bfz}-{\bfa}(s))+{\bfa}(t),$$
 where 
$R(s,t)$ is the rotation operator with the angle between ${\bfa}(s)$ and 
${\bfa}(t)$. Then
 $$ \begin{array}{l}
  \displaystyle{|\tilde {\bfz}-{\bfa}(t)|= \left(1-\frac 1{\gamma_*}
   |\gamma({\boldsymbol\zeta}(s))-\gamma({\boldsymbol\zeta}(t))
   |\right)|{\bfz}-{\bfa}(s)| }\\[0.3cm]
  \displaystyle{~~~~~~~~~~~~~\le \left(1-\frac 1{\gamma_*}
   |\gamma({\boldsymbol\zeta}(s))-\gamma({\boldsymbol\zeta}(t))
   |\right)\gamma({\boldsymbol\zeta}(s))}\\[0.3cm]
  \displaystyle{~~~~~~~~~~~~~
    = \gamma({\boldsymbol\zeta}(s))-\frac{\gamma({\boldsymbol\zeta}(s))}
     {\gamma_*}
  |\gamma({\boldsymbol\zeta}(s))-\gamma({\boldsymbol\zeta}(t))|\le 
   \gamma({\boldsymbol\zeta}(t)),}
 \end{array} $$
and 
  \begin{eqnarray*}
   (\tilde{\bfz}-{\bfa}(t))\cdot {\bfa}(t)&=&\left(1-\frac 1{\gamma_*}|
  \gamma({\boldsymbol\zeta}(s))-\gamma({\boldsymbol\zeta}(t))|\right)
  R(s,t)({\bfz}-{\bfa}(s))\cdot{\bfa}(t)\\
    &=&\left(1-\frac 1{\gamma_*}|
  \gamma({\boldsymbol\zeta}(s))-\gamma({\boldsymbol\zeta}(t))|\right)
  R(s,t)({\bfz}-{\bfa}(s))\cdot R(s,t)^{-1}{\bfa}(s)\\
    &=&\left(1-\frac 1{\gamma_*}|
  \gamma({\boldsymbol\zeta}(s))-\gamma({\boldsymbol\zeta}(t))|\right)
   ({\bfz}-{\bfa}(s))\cdot {\bfa}(s)=0.
  \end{eqnarray*}
Hence $\tilde{\bfz} \in K(\theta;t)$.
Besides, since $|R(s,t)-I|_{B({\bf R}^2)} \le C_{{\bfa}}|{\bfa}(t)-{\bfa}(s)|$
for all $s,t \in [0,T]$ and 
for some constant $C_{\bfa}>0$ depending only on ${\bfa}$ (cf. \eqref{(6.1)}-\eqref{(6.3)}),
we see that
  \begin{eqnarray*}
  |\tilde {\bfz}-{\bfz}| &\le&\left| \left(1-\frac 1{\gamma_*}|
  \gamma({\boldsymbol\zeta}(s))-\gamma({\boldsymbol\zeta}(t))|\right)
  R(s,t)({\bfz}-{\bfa}(s))+{\bfa}(t)-{\bfz}\right|\\ 
   &\le& |R(s,t)({\bfz}-{\bfa}(s)) -({\bfz}-{\bfa}(s))|+|{\bfa}(t)-{\bfa}(s)|\\
   && ~~~+\frac 1{\gamma_*}\left|\gamma({\boldsymbol\zeta}(s))
     -\gamma({\boldsymbol\zeta}(t)) \right| |{\bfz}-{\bfa}(s)|\\
   &\le& {\rm const.}\int_s^t \left\{|{\bfa}'(\tau)|+\left|\frac d{d\tau}
    \gamma({\boldsymbol\zeta}(\tau)\right|\right\}d\tau.
   \end{eqnarray*}
Similarly, for any $p \geq 2$,
 $$ \frac 1p|\tilde {\bfz}|^p-\frac 1p|{\bfz}|^p \le {\rm const.}
    |\tilde{\bfz}-{\bfz}| \le 
    {\rm const.}\int_s^t\left\{|{\bfa}'(\tau)|+\left|\frac d{d\tau}
    \gamma({\boldsymbol\zeta}(\tau))\right| \right\}d\tau. $$
For any $s,~t \in [0,T]$ we get the same type of estimates as above by 
repeating the above procedure at most $N$-times. Thus we can see that 
$\Theta_1 \subset \Theta_S$, since 
$\gamma({\boldsymbol\zeta})$
and ${\bfa}$ are of $C^2$. Hence $\Theta_1 \subset \Theta_W\subset
\Theta$. Since $\Theta_1$ is dense in $\Theta$, we conclude that
$\Theta_W=\Theta$ and $(\Theta_S)$ holds.

\medskip

\noindent
({\it Verification of (A1)-(A3)}) For any parameters $\theta:=[{\bfa},\gamma,
{\boldsymbol\zeta}],
~\bar\theta:=[\bar{\bfa},\bar\gamma,\bar{\boldsymbol\zeta}] \in \Theta$, 
we have by our assumption that
${\bfa},~\bar{\bfa} \in W^{1,p}(0,T;{\bf R}^2)$, so that expression 
\eqref{(6.1)}-\eqref{(6.3)}
of $R_{\theta\bar\theta}$ shows that $R_{\theta\bar\theta} \in 
W^{1,p}(0,T; B({\bf R}^2))$ and (A1) holds. 
Since $\bar{\bfa}(t)=
R_{\theta\bar\theta}(t){\bfa}(t)$, 
$\bar {\bfz}:=(1-\varepsilon)R_{\theta\bar\theta,\varepsilon}{\bfz}
   +\varepsilon \bar{\bfa}(t),~0<\varepsilon<1,$ is written as
$$ \bar{\bfz}-\bar{\bfa}(t)=(1-\varepsilon)R_{\theta\bar\theta}
   ({\bfz}-{\bfa}(t)).$$
Hence, 
 $$ (\bar{\bfz}-\bar{\bfa}(t))\cdot \bar{\bfa}(t)
   =(1-\varepsilon)R_{\theta\bar\theta}(t)({\bfz}-{\bfa}(t))\cdot
    R_{\theta\bar\theta}(t){\bfz}
    =(1-\varepsilon)({\bfz}-{\bfa}(t))\cdot {\bfz}=0.$$
and if $d_{\Theta}(\theta,\bar\theta) <\varepsilon\gamma_*$, then
 $$|\bar{\bfz}-\bar{\bfa}(t)|\le (1-\varepsilon)\gamma({\boldsymbol\zeta}(t))
   \le \gamma({\boldsymbol\zeta}(t))-\varepsilon \gamma_* \le 
       \bar\gamma(\bar{\boldsymbol\zeta}(t)).$$
and $\bar{\bfz} \in K(\bar\theta;t)$. By \eqref{(6.4)} and \eqref{(6.5)}, we see that (A2) and 
(A3) with 
$\delta_\varepsilon=\varepsilon\gamma_*$ are fulfilled.

By virtue of Theorem \ref{th:3.1}, the time-derivative operator 
$L_{{\bfu}_0}(\theta;\cdot)$ is defined as a maximal monotone mapping 
from $ L^p(0,T;{\bf R}^2)$ into $L^{p'}(0,T;{\bf R}^2)$ associated with 
$\{K(\theta;t)\}$ given by \eqref{(6.5)} and an initial
datum ${\bfu}_0 \in \overline {K(\theta;0)}$. By our assumption, 
$\cup_{\theta\in \Theta, t\in [0,T]} K(\theta;t)$ is bounded in ${\bf R}^2$,
namely there is a closed ball $B_{k_0}$ around the origin of ${\bf R}^2$
and with radius $k_0>0$ such that
    $$\bigcup_{\theta\in \Theta, t\in [0,T]}  K(\theta;t) \subset B_{k_0}.$$

\medskip

\noindent
({\it Feedback system} $\Lambda_{{\bfu}_0}$) 
We define a feedback system $\Lambda_{{\bfu}_0}$ as follows.
Let ${\bfG}:={\bfG}(t,{\bfw},{\boldsymbol \zeta})$ be a globally bounded and 
continuous vector field from $[0,T]\times {\bf R}^2\times {\bf R}^2$ into 
${\bf R}^2$ such that
 $$ |{\bfG}(t,{\bfw},{\boldsymbol \zeta}) -{\bfG}(t,\bar{\bfw},
   \bar{\boldsymbol \zeta})| \le C_{{\bfG}}(|{\bfw}-\bar{\bfw}|+ 
   |{\boldsymbol \zeta} -\bar{\boldsymbol \zeta}|), ~~\forall
   {\bfw},~\bar{\bfw},~{\boldsymbol \zeta},
   ~\bar{\boldsymbol \zeta}\in {\bf R}^2,~\forall t\in [0,T],$$
where $C_{\bfG}$ is a positive constant. Let ${\bfY}$ be a closed convex
set in ${\bf R}^2$ such that $0 \notin {\bfY}$. 
Now, given ${\boldsymbol \zeta}_0 \in {\bfY}$ and
${\bfu} \in {\cal U}:=\{{\bfw} \in L^p(0,T;{\bf R}^2)~|~ {\bfw}(t) \in B_{k_0},
{\rm a.e.}~t \in (0,T)\}$, consider the evolution inclusion
\begin{equation} {\boldsymbol\zeta}'(t)+\partial I_{\bfY}({\boldsymbol\zeta}(t))
    \ni {\bfG}(t,\int_0^t{\bfu}(\tau)d\tau, {\boldsymbol\zeta}(t)),~
    {\rm a.e.~}t\in (0,T),~~ {\boldsymbol\zeta}(0)={\boldsymbol\zeta}_0.
  \label{(6.6)}
    \end{equation}
By the general theory of nonlinear evolution equations (cf. [4]) \eqref{(6.6)} has a 
unique solution in $W^{1,2}(0,T;{\bf R}^2)$ and for
simplicity this solution is denoted by ${\boldsymbol\zeta}({\bfu})$ or 
${\boldsymbol\zeta}({\bfu};t)$.

\begin{lem}\label{lem:6.1}
Let $\{{\bfu}_n\}$ be a sequence in ${\cal U}$ such that
${\bfu}_n \to {\bfu}$ weakly in $L^p(0,T;{\bf R}^2)$ (as $n \to \infty$), then
the sequence ${\boldsymbol\zeta}({\bfu_n})$ converges to
the solution ${\boldsymbol\zeta}({\bfu})$ of \eqref{(6.6)} in the sense that
\begin{equation} {\boldsymbol\zeta}_n :={\boldsymbol\zeta}({\bfu}_n) \to
   {\boldsymbol\zeta}:={\boldsymbol\zeta}({\bfu})
    ~{\it in~}C([0,T];{\bf R}^2),~~
  {\boldsymbol\zeta}'_n \to {\boldsymbol\zeta}'~{\it in~}
   L^p(0,T;{\bf R}^2). \label{(6.7)}
    \end{equation}  
\end{lem}

\noindent
{\bf Proof.} Multiplying the both sides of \eqref{(6.6)} for ${\boldsymbol\zeta}_n$
by ${\boldsymbol\zeta}'_n$, we get
  $$|{\boldsymbol\zeta}'_n(t)|^2
  + \frac d{dt}I_{\bfY}({\boldsymbol\zeta}_n(t))
    ={\bfG}_n(t)\cdot {\boldsymbol\zeta}'_n(t),~~{\rm a.e.~}t \in (0,T),$$
where ${\bfG}_n(t):={\bfG}(t,\int_0^t {\bfu}_n(\tau)d\tau,
{\boldsymbol\zeta}_n(t))$.
Since $I_{\bfY}({\boldsymbol\zeta}_n) \equiv 0$ on $[0,T]$, we see that
$|{\boldsymbol\zeta}'_n(t)| \le |{\bfG}_n(t)|$ and hence 
$\{{\boldsymbol\zeta}_n\}$ is 
bounded in
$W^{1,\infty}(0,T;{\bf R}^2)$. Therefore there exists a subsequence 
$\{{\boldsymbol\zeta}_{n_k}\}$ of $\{{\boldsymbol\zeta}_n\}$ such that 
 $$ {\boldsymbol\zeta}_{n_k} \to {\boldsymbol\zeta}~
   {\rm in~}C([0,T];{\bf R}^2),~~
  {\boldsymbol\zeta}'_n \to {\boldsymbol\zeta}'
   ~{\rm weakly~in~}L^p(0,T;{\bf R}^2).$$
for some function ${\boldsymbol\zeta}\in W^{1,\infty}(0,T;{\bf R}^2)$. Since
${\bfG}_n \to {\bfG}:={\bfG}(t,\int_0^t{\bfu}d\tau,{\boldsymbol\zeta})$ in 
$C([0,T];{\bf R}^2)$, it follows that ${\boldsymbol\zeta}_{n_k}$ converges in
$C([0,T];{\bf R}^2)$ to the
solution of \eqref{(6.6)} which is nothing but ${\boldsymbol\zeta}$ and consequentlyy
${\boldsymbol\zeta}_n \to {\boldsymbol\zeta}$ in $C([0,T];{\bf R}^2)$ without
extracting any subsequence from $\{{\boldsymbol\zeta}_n\}$ by the uniquenss of
solution to \eqref{(6.6)}. Besides we 
have for any $s,~t \in [0,T],~s\le t$, 
  \begin{eqnarray*} 
   \limsup_{n\to \infty}\int_s^t|{\boldsymbol\zeta}'_n|^2d\tau
    &=& \lim_{n\to \infty}\int_s^t G_n\cdot {\boldsymbol\zeta}'_n d\tau\\
      &=&\int_0^T G\cdot {\boldsymbol\zeta}' dt
    =\int_s^t|{\boldsymbol\zeta}'|^2 d\tau,
  \end{eqnarray*}
so that ${\boldsymbol\zeta}'_n \to {\boldsymbol\zeta}'$ in 
$L^2(0,T;{\bf R}^2)$ (hence in $L^p(0,T;{\bf R}^2)$) as well as 
$|{\boldsymbol\zeta}'(t)|\le |{\bfG}(t)|$
for a.e. $t\in (0,T)$. 
Thus \eqref{(6.7)} has been obtained. \hfill $\Box$ \vspace{0.5cm}

Now, we define a mapping 
${\bfa}(\cdot)$ by putting 
  $${\bfa}({\bfu};t) :=
    \frac {{\boldsymbol \zeta}({\bfu};(t))}
 {|{\boldsymbol \zeta}({\bfu}; t)|}$$
for any ${\bfu}$ in ${\cal U}$. It is easy to see that
 $${\bfa}({\bfu}) \in W^{1,p}(0,T;{\bf R}^2),~~
 |{\bfa}({\bfu};(t)|=1,~\forall t \in [0,T],~{\bfu} \in {\cal U}.$$
Now our feedback system $\Lambda_{{\bfu}_0}: {\cal U} \to
\Theta$ is given by 
 $$ \Lambda_{{\bfu}_0}{\bfu}=[{\bfa}({\bfu}),\gamma,
 {\boldsymbol\zeta}({\bfu})],~~\forall {\bfu}\in {\cal U}.$$

We observe from Lemma \ref{lem:6.1} that $\Lambda_{{\bfu}_0}$ is a compact mapping from 
${\cal U}$ into $\Theta$; in fact, if ${\bfu}_n \in {\cal U}$ and 
${\bfu}_n \to 
{\bfu}$ weakly in $L^p(0,T;{\bf R}^2)$ (as $n \to \infty$), then 
${\boldsymbol\zeta}({\bfu}_n) \to {\boldsymbol\zeta}({\bfu})$ in 
$W^{1,p}(0,T;{\bf R}^2)$ and $\theta_n:=
[{\bfa}({\bfu}_n),\gamma,{\boldsymbol\zeta}({\bfu}_n)] \to \theta:=
[{\bfa}({\bfu}),\gamma,{\boldsymbol\zeta}({\bfu})]$ in $\Theta$.\vspace{0.5cm}

In the above set-up, given an initial datum ${\bfu}_0$, satisfying that
 $$ {\bfa}_0:=\frac {{\boldsymbol\zeta}_0}{|{\boldsymbol\zeta}_0|},
  ~{\bfa}_0\cdot ({\bfu}_0-{\bfa}_0)=0,~
    |{\bfu}_0 -{\bfa}_0|\le \gamma({\boldsymbol\zeta}_0),$$
and ${\bff} \in L^{p'}(0,T;{\bf R}^2)$, we formulate 
a quasi-variational ordinary differential inequality: 
\begin{equation}  L_{{\bfu}_0}(\theta;{\bfu})\ni {\bff}
         ~~{\rm in~}~L^{p'}(0,T;{\bf R}^2),~~
     \theta=\Lambda_{{\bfu}_0}{\bfu} \label{(6.8)}
    \end{equation}
which is written as
 \begin{eqnarray*}
 &&  {\boldsymbol\zeta}'(t)+\partial I_{\bfY}({\boldsymbol\zeta}(t))
   \ni {\bfG}(t, \int_0^t{\bfu}(\tau)d\tau,{\boldsymbol\zeta}(t)),~t\in (0,T),
     ~~{\boldsymbol\zeta}(0)={\boldsymbol\zeta}_0,\\
&&  \int_0^T ({\boldsymbol\eta}'(t)-{\bff}(t)) \cdot ({\bfu}()
     -{\boldsymbol\eta}(t))dt
    \le \frac 12|{\bfu}_0-{\boldsymbol\eta}|^2,~\forall {\boldsymbol\eta}\in 
    {\cal K}_0(\theta).
 \end{eqnarray*}

The existence of a solution $\{\theta,{\bfu}\}$ of \eqref{(6.8)} is not covered
by Theorem \ref{th:5.1}, because condition \eqref{(5.3)} is not fulfilled. Hence the compactness
result mentioned in Lemma \ref{lem:5.1} is not obtained. However the 
existence of a solution of \eqref{(6.8)} is 
directly proved by the fixed point argument. In fact, let $S$ be the mapping
which assigns to each function ${\bfu}\in {\cal U}$ the solution 
$\bar{\bfu}$ of
 $$\theta= \Lambda_{{\bfu}_0}{\bfu},~~L_{{\bfu}_0}(\theta;\bar{\bfu})
    \ni {\bff}.$$
Then it follows from the above observations that $S$ is compact mapping from
${\cal U}$ into itself in the topology of $L^p(0,T;{\bf R}^2)$, so that 
$S$ has at
least one fixed point, ${\bfu}=S{\bfu}$ in ${\cal U}$. Clearly, the pair 
$\{\theta,{\bfu}\}$ with 
$\theta:=[{\bfa}({\bfu}),\gamma,{\boldsymbol\zeta}({\bfu})]$
is a solution of \eqref{(6.8)}. 

\begin{req}\label{rem:6.1}
For simplicity we treated above in ${\bf R}^2$ . 
But similar problems are formulated and solved in the 3d
space or more generally infinite dimensional spaces, too, although the 
computation is more technical for the verification of assumptions (A1)-(A3)
and $(\Theta_S)$.
\end{req}

\subsection{Quasi-variational partial differential inequality}
\label{sec:superconductor}

In this application we treat a model arising in superconductivity.
Let $\Omega$ be a smooth bounded domain in ${\bf R}^N,~1\le N <\infty$, and
  $$\Gamma:=\partial \Omega,~~\Sigma:=\Gamma\times (0,T),~~
     Q:=\Omega \times (0,T).$$
We put $ V:=W_0^{1,p}(\Omega),~H:=L^2(\Omega),
~W:=W_0^{2,q}(\Omega)$ with $\max\{p,N\} <q <\infty$, and hence 
$W\subset V \subset H \subset V^*\subset W^*$ with dense and compact 
embeddings. We suppose that
\begin{itemize}
\item $a(x,t,v)$ is a function on $Q \times {\bf R}$, satisfying the
Carath\'eodory condition, namely for a.e. $(x,t) \in Q$, the function $v \to
a(x,t,v)$ is continuous and for all $v \in {\bf R}$ the function $(x,t) \to
a(x,t,v)$ is measurable on $Q$. We assume that for some positive constants 
$a_*,~a^*$
  $$ a_* \le a(x,t,v) \le a^*,~~{\rm a.e.~}(x,t)\in Q,~\forall v \in {\bf R}.
    $$

\item $\gamma(\cdot)\in C_b({\bf R})$, namely $\gamma$ is bounded and 
continuous on ${\bf R}$ such that $\lim_{r \to \pm \infty}\gamma(r)$ exists.
Furthermore, suppose that, for a positive constant $\varepsilon_0$,
$$  \gamma(\zeta) \geq \varepsilon_0,~~\forall \zeta \in {\bf R}.$$

\item $h(x,t,u)$ is a Lipschitz continuous function on 
$\overline{Q}\times {\bf R}$.
\end{itemize}

Next, we put
 $$ \Theta:=\{\theta:=[\gamma, \zeta]~|~\gamma\in C_b({\bf R}),
  \gamma \geq \varepsilon_0~{\rm on~}{\bf R},~\zeta \in C(\overline Q)\},
  $$ 
and define a semimonotone mapping $A(v;u)$ and a family $\{K(\theta;t)\}$ 
of closed and convex subsets of $V$ by:
 $$ A(v;u):= -{\rm div}\hspace{0.05cm}(a(x,t,v)|\nabla u|^{p-2}\nabla u)
              ~~{\rm in~} L^{p'}(0,T;V^*), $$
     $$~~~~~~~~~~\forall v \in {\cal V}_A:=L^p(0,T;V),~\forall u \in L^p(0,T;V),$$
 $$ K(\theta;t):=\{z \in V~|~|\nabla z(x)| \le \gamma(\zeta(x,t))~{\rm a.e.~}
    x \in \Omega\}, $$
 $$ ~~~~~~~~~~~\forall t \in [0,T],~\forall \theta:=[\gamma,\zeta] \in\Theta.
         $$ 
Furthermore, for a given $u_0 \in H$, a feedback system $\Lambda_{u_0}$ is
defined by: $\theta:=[\gamma,\zeta]=\Lambda_{u_0}u$ if and only if 
$u \in L^p(0,T;V)$ and
$\zeta$ is a unique solution of the following heat equation
\begin{equation}  \left\{ \begin{array}{l}
  \displaystyle{\zeta_t -\Delta \zeta= h(x,t, u)~~{\rm in~}Q,}\\[0.3cm]
  \displaystyle{ \frac {\partial\zeta}{\partial n}+n_0 \zeta=0~~{\rm on~}
     \Sigma,~~\zeta(\cdot,0)=\zeta_0~{\rm in~}\Omega,}
   \end{array} \right. \label{(6.80)}
    \end{equation}
where $\zeta_0 \in H^2(\Omega)$, satisfying 
$\frac {\partial\zeta_0}{\partial n}+n_0 \zeta_0=0~{\rm on~}\Gamma$ and
$u_0 \in K(\theta_0;0):=\{z \in V~|~|\nabla z|\le \gamma(\zeta_0)~
{\rm a.e.~on~}\Omega\}$.\vspace{0.5cm}

We know (cf. [7; Appendix]) that
problem \eqref{(6.80)} admits a unique solution 
$\zeta \in W^{1,2}(0,T;H^1(\Omega))\cap L^\infty(0,T; H^2(\Omega))$, 
which is compactly embedded in 
$C(\overline {Q})$. \vspace{0.5cm}

It is verified that $(\Theta_S)$ holds under \eqref{(5.2)}-\eqref{(5.3)}, 
$\Theta\cap (C^2({\bf R})\times C^2(\overline Q)) \subset \Theta_S$, 
$\Theta_W=\Theta$ and conditions 
(A1)-(A3) are satisfied by $R_{\theta\bar\theta}=I$, 
$\sigma_{\theta\bar\theta,\varepsilon}=0$ and 
${\cal F}_{\theta\bar\theta,\varepsilon}=(1+c_0(\varepsilon))I$;
see [7, 20] for the verification of these facts. 
Also, it is easy to check conditions (B0)-(B3) from the definition of $A(v;u)$
 as well as $\Lambda_{u_0}$ satisfies $(\Lambda1)-(\Lambda2)$. 
Therefore, by virtue of Theorem \ref{th:5.1}, for a given 
$f \in L^{p'}(0,T;V^*)$, the quasi-variational inequality
  $$ L_{u_0}(\theta;u)+A(u;u) \ni f~~{\rm in~}L^{p'}(0,T;V^*),
     ~~\theta=\Lambda_{u_0}u, $$
has at least one solution $\{\theta,u\}$, which gives a solution of the 
following system:
 \begin{eqnarray*}
 &&  \zeta_t -\Delta \zeta= h(x,t, u)~~{\rm in~}Q,~~
  \frac {\partial\zeta}{\partial n}+n_0 \zeta=0~{\rm on~}
     \Sigma,~~\zeta(\cdot,0)=\zeta_0~{\rm in~}\Omega, \\
 && u \in C([0,T]; H),~~u(t) \in K(\theta;t),~\forall t \in [0,T],\\
 && \int_0^T \langle \eta',u-\eta\rangle dt 
    +\int_Q a(x,t,u)|\nabla u|^{p-2}\nabla u\cdot \nabla(u-\eta)dxdt\\
  && ~~~~~
    \le \int_0^T\langle f,u-\eta\rangle dt +\frac 12|u_0-\eta(0)|^2_H,
   ~~\forall \eta \in {\cal K}_0(\theta).
 \end{eqnarray*}

\begin{req}
 The above quasi-variational inequality arises from 
simplified models of\\
 Type II superconductivity; we refer to [1, 2, 19, 27] for related works.
\end{req}

\begin{req}
The similar approach is possible to quasi-variational Navier-Stokes problems. 
In the general case the obstacle function $\gamma(\zeta)$ is required to 
have three phases, $\gamma(\zeta)=0$, $0< \gamma(\zeta)<\infty$ and 
$\gamma(\zeta)=\infty$ which are respectively the solid, mussy and liquid
parts in the fluid. Therefore, its mathematical treatment would be much more 
complicated
(cf. [9, 11, 12]).
\end{req}


\begin{center}
{\bf Acnowledgment}
\end{center}

We thank the Cardinal Stefan Wyszy\'nski University in Warsaw 
for supporting our collaboration.
Many thanks to Jakub Zieli\'nski and Krzysztof Nowi\'nski from ICM, 
University of Warsaw, for important discussions.

\begin{center}
{\bf References}
\end{center}

\begin{enumerate}
\item G. Akagi, Convergence of functionals and its applcations to parabolic
equations, Abst. Appl. Anal., 11(2004), 907-933.
\item A. Azevedo and L. Santos, A diffusion problem with gradient 
constraint depending on the temperature, Adv. Math. Sci. Appl., 20(2010),
151-166.
\item A. Bensoussan and J. L. Lions, Contr\^ole impulsionnel et controle
contitinu m\'ethode des in\'equations quasi variationnelles non lin\'eaires,
C. R. Acad. Sci. Paris, 278(1974), 675-679. 

\item H. Br\'ezis, {\it Op\'erateurs maximaux monotones et semi-groupes de 
contractions dans les espaces de Hilbert}, Math. 
Studies 5, North-Holland, Amsterdam, 1973.

\item M. Cyrot and D. Pavuna,  Introduction to superconductivity and high-Tc materials, World Scientific Publishing Company, 1992.

\item F. E. Browder, Nonlinear maximal monotone operators in Banach spaces,
Math. Ann., 175(1968), 89-113.

\item T. Fukao and N. Kenmochi, Parabolic variational inequalities with
weakly time-dependent constraints, Adv. Math. Sci. Appl., 23(2013), 365-395.
\item T. Fukao and N. Kenmochi, Quasi-variational inequality approach to
heat convection problems with temperature dependent velocity constraint,
Discrete Contin. Dyn. Syst., 6(2015), 2523-2538.
\item M. Gokieli, N. Kenmochi and M. Niezg\'odka, 
Variational inequalities of Navier--Stokes type with time dependent 
constraints, J. Math. Anal. Appl., 449(2017), 1229-1247
\item M. Gokieli, N. Kenmochi and M. Niezg\'odka, A new compactness theorem 
for variational inequalities of parabolic type, Houston J. Math. 44(2018), 
319-350.
\item M. Gokieli, N. Kenmochi and M. Niezg\'odka, Mathematical modelling of 
biofilm development, Nonlinear Anal. Real World Anal., 42(2018), 422-447.
\item M. Gokieli, N. Kenmochi and M. Niezg\'odka, Variational and 
quasi-variational inequalities of Navier-Stokes type with velocity constraints,
Adv. Math. Sci. Appl., 27(2018), 359-402.
\item M. Hinterm\"uller and C. Rautenberg, Parabolic variational inequalities
with gradient-type constraints, SIAM J. Optim., 23(2013), 2090-2123.
\item M. Hinterm\"uller and C. Rautenberg, On the uniqueness and numerical
approximation of solutions to certain parabolic quasi-variational inequalities, Port. Math., 74(2017), 1-35.

\item J.-L. Joly and U. Mosco, A propots de l'existence et de la r\'gularit\'e
des solutions de certaines quasi-variaionnelles, J. Funct. Anal., 34(1979),
107-137.

\item A. Kadoya, Y. Murase and N. Kenmochi, A class of nonlinear parabolic
systems with environmental constraints, Adv. Math. Sci. Appl., 20(2010), 
281-313.
\item R. Kano, N. Kenmochi and Y. Murase, Nonlinear evolution equations
generated by subdifferentials with nonlocal constraints, Banach Center
Publication, 86(2009), 175-194.
\item N. Kenmochi, Solvability of nonlinear evolution equations with
time-dependent constraits and applications, Bull. Fac. Education, Chiba
Univ., 30(1981), 1-87.
\item N. Kenmochi, Parabolic quasi-variational diffusion problems with
gradient constraints, Discrete. Contin. Dyn. Syst. 6(2013), 423-438.
\item N. Kenmochi and M. Niezg\'odka, Weak solvability for parabolic 
variational inclusions and applications to quasi-variational inequalities,
Adv. Math. Sci. Appl., 25(2016), 62-97.

\item M. Kunze and J. F. Rodrigues, An elliptic quasi-variational inequality
with gradient constraints and some of applications, Math. Methods Appl. Sci.,
23(2000), 897-908.
\item J. Leray and J. L. Lions, Quelques r\'esultats de Visik sur les 
probl\`emes elliptique non lin\'eaires par les m\'ethodes de Minty-Browder,
Bull. Soc. Math. France, 93(1965), 97-107.
\item F. Mignot and J.-P. Puel, In\'equations d'\'evolution parab\'equations 
avec convexes d\'ependant du temp. Applications aux in\'equations quasi
variationelles d'\'evolution, Arch. Rattional Mech. Anal., 64(1977), 493-519.
\item U. Mosco, On the continuity of the Young-Fenchel transform, J. Math. 
Anal. Appl., 35(1971), 518-535.
\item Y. Murase, Abstract quasi-variational inequalities of elliptic type and
applications, Banach Center Publication, 86(2009), 175-194.235-246.
\item L. Prigozhin, On the Bean critical-state model in superconductivity,
European J. Appl. Math., 7(1996), 237-247.
\item J.-F. Rodrigues and L. Santos, A parabolic quasi-variational inequality
arising in a superconductivity model, Ann. Scuola Norm. Sup. Pisa, Sci.,
29(2000), 153-169.
\item U. Stefanelli, Nonlocal quasi-variational evolution equations, J. 
Differential Equations, 229(2006), 204-228.
\item L. Tartar, In\'equations quasi variationnelles abstraites, C. R. Acad. 
Sci. Pris, 278(1974), 1193-1196.
\item Y. Yamada, On evolution equations generated by subdifferential operators,
J. Fac. Sci. Univ. Tokyo Sect. IA., 23(1976), 491-515.
\end{enumerate}

\end{document}